\magnification 1200
\input amstex
\documentstyle{amsppt}
\vsize=7.85in
\topmatter

\rightheadtext{}
\leftheadtext{M. S. Baouendi,  P.
Ebenfelt, and L. P. Rothschild}
\title Convergence and finite determination of formal
CR mappings\endtitle
\author M. S. Baouendi,  P. Ebenfelt, and Linda
Preiss Rothschild
\endauthor
\address Department of Mathematics, 0112, University of
California at San Diego, La Jolla, CA
92093-0112\endaddress
\email sbaouendi\@ucsd.edu,
lrothschild\@ucsd.edu\endemail
\address Department of Mathematics, Royal
Institute of Technology, S-100 44 Stockholm,
Sweden\endaddress
\email ebenfelt\@math.kth.se\endemail
\thanks The first and the third authors are
partially supported by National Science
Foundation grant DMS 98-01258. The second author
is partially supported by a grant from the Swedish Natural
Science Research Council.\endthanks

\subjclass 32H02\endsubjclass
\loadeufm

\def\det{{\text{\rm det}}}

\define \p{\eufm p}
\define \q{\eufm q}
\define \m{\eufm m}
\define \im{\text{\rm Im }}
\define \re{\text{\rm Re }}


\define \bR{\Bbb R}
\define \bC{\Bbb C}



\def\rk{\text{\rm rk}}

\define \dbl{[\![}
\define \dbr{]\!]}
\define \hO{\bC\dbl x\dbr}

\abstract It is shown that a formal mapping between two
real-analytic hypersurfaces in complex space is
convergent provided that neither hypersurface contains a
nontrivial holomorphic variety. For higher codimensional
generic submanifolds, convergence is proved e.g.\ under the
assumption that the source is of finite type, the target does
not contain a nontrivial holomorphic variety, and the mapping
is finite. Finite
determination (by jets of a predetermined order) of
formal mappings between smooth generic submanifolds is also
established.
\endabstract

\endtopmatter

\document

\head \S 0. Introduction\endhead

In this paper, we study the convergence and finite
determination of formal holomorphic mappings of
$(\bC^N,0)$ taking one real submanifold into
another. By a formal (holomorphic) mapping
$H\:(\bC^N,0)\to(\bC^N,0)$, we mean an $N$-vector
$H=(H_1,\ldots, H_N)$, where each $H_j$ is a formal power
series in $N$ indeterminates with no constant term. If
$M$ and
$M'$ are real smooth submanifolds through $0$ in $\bC^N$
defined near the origin by
$\rho(Z,\bar Z)=0$ and $\rho'(Z,\bar Z)=0$ respectively,
where $\rho$ and $\rho'$ are vector valued smooth defining
functions, then we say that a formal mapping
$H\:(\bC^N,0)\to(\bC^N,0)$ sends $M$ into $M'$ if the
vector valued  power series $\rho'(H(Z),\overline{H(Z)})$ is
a (matrix) multiple of $\rho(Z,\bar Z)$. For real-analytic 
hypersurfaces, we shall prove the following.

\proclaim{Theorem 1} Let $M$ and $M'$ be real-analytic
hypersurfaces through the origin in $\bC^N$, $N\geq 2$. Assume
that neither
$M$ nor $M'$ contains a nontrivial holomorphic
subvariety through $0$. Then any formal mapping
$H\:(\bC^N,0)\to (\bC^N,0)$ sending $M$ into $M'$ is
convergent.   
\endproclaim 

The condition that $M'$ above does not
contain a nontrivial holomorphic subvariety is necessary (see
Remark 2.3). As a corollary we obtain the following 
characterization. 

\proclaim{Corollary 2} Let $M$ be a real-analytic
hypersurface through $0$ in $\bC^N$ with $N\geq 2$. Then, $M$
does not contain a nontrivial holomorphic subvariety through
$0$ if and only if every formal mapping $H\:(\bC^N,0)\to
(\bC^N,0)$ sending $M$ into itself is convergent.\endproclaim

As a consequence of the geometric properties of mappings
used in the proof of Theorem 1 above, we also obtain a new
reflection principle for CR mappings between real-analytic
hypersurfaces (see Theorem 2.4). 

We state now some results for 
generic submanifolds of higher codimension in $\bC^N$. A
formal  mapping $H\:(\bC^N,0)\to (\bC^N,0)$ is called {\it
finite} if the ideal generated by its components
$H_1,\ldots, H_N$ has finite codimension in the ring of
formal power series $\bC\dbl Z\dbr$. 
Recall that a
smooth generic submanifold $M\subset\bC^N$ is of {\it finite
type} (in the sense of Kohn and Bloom--Graham) at $p\in M$ if
the Lie algebra generated by its CR vector fields and their
complex conjugates has full dimension at $p$ (see \S 1.1). For
a real-analytic hypersurface in $\bC^2$ the notion of finite
type is equivalent to that of not containing a nontrivial
variety. Our first result for higher codimensional generic
submanifolds is the following.

\proclaim{Theorem 3} Let $M$ and
$M'$ be real-analytic generic submanifolds of the same
dimension through the origin in $\bC^N$. Assume
that
$M$ is of finite type at $0$, and that $M'$ does not contain
any germ of a nontrivial holomorphic subvariety through $0$.
Then any  formal finite mapping $H\:(\bC^N,0)\to
(\bC^N,0)$ sending $M$ into $M'$ is convergent. 
\endproclaim

In fact, we prove a more general result (see Theorem 2.1) in
which $M$ and $M'$ need not be in the same complex space nor of
the same dimension. Moreover, the condition that
$M'$ does not contain any nontrivial variety is replaced by the
weaker assumption that $M'$ is essentially finite at $0$, and
the condition of finiteness of the mapping $H$ is replaced by
injectitivity of its Segre homomorphism. (See \S 1.1 for
relevant definitions.) 

  Recall that two real submanifolds $M$ and
$M'$ of the same dimension through the origin in
$\bC^N$ are said to be {\it biholomorphically
equivalent} at $0$ if there exists a biholomorphism
$H\:(\bC^N,0)\to (\bC^N,0)$, defined near $0$, sending $M$
into $M'$. Similarly, we shall say that $M$ and $M'$ are {\it
formally equivalent} (at $0$) if there exists a formal
invertible mapping  sending $M$ into $M'$. 

\proclaim{Theorem 4} Let $M$ and
$M'$ be real-analytic generic submanifolds
through the origin in $\bC^N$. Assume
that
$M$ is of finite type at $0$ and does not contain any germ of a
nontrivial holomorphic subvariety through $0$. Then $M$ and
$M'$ are formally equivalent at $0$ if and only if they are
biholomorphically equivalent at $0$.
\endproclaim

In the following result, we use
the notation $j^k_p f$ for the jet of order $k$ at $p\in M$ of
a smooth mapping $f\: M\to M'$, where 
$M$ and $M'$ are smooth manifolds. We shall say that a
smooth generic submanifold
$M$ does not contain a (nontrivial) formal variety through $0$
if there is no (nontrivial) formal holomorphic curve
$\gamma\:(\bC,0)\to (\bC^N,0)$ such that
$\rho(\gamma(t),\overline{\gamma(t)})\sim 0$, where the
composition is in the sense of formal power series and $\sim$
means equality for power series. Here,
$\rho(Z,\bar Z)$ is a vector valued smooth defining function
for $M$ near $0$ identified with its Taylor series at the
origin. 

\proclaim{Theorem 5} Let $M$ be a smooth generic
submanifold through the origin in $\bC^N$. Assume that $M$ is 
of finite type at $0$ and does not contain a nontrivial
formal variety through
$0$. Then there exists a positive integer $K$ with the
following property. If $M'$ is a smooth generic
submanifold in $\bC^N$ and  $H^1,H^2\:M\to
M'$ are smooth
local CR diffeomorphisms, defined near
$0$, such that $j^K_0(H^1)=j^K_0(H^2)$, then
$j^l_0(H^1)=j^l_0(H^2)$ for all positive integers $l$.  
\endproclaim 

Theorem 5 is a consequence of the more general result
Theorem 2.5, in which $M$ and $M'$ can be of different
dimensions and be contained in different complex spaces. As in
the case of Theorem 3, we may replace the condition that
$M$ does not contain a nontrivial formal variety through
$0$ by the assumption that $M$ is
essentially finite at $0$ and the condition that $H^1$,
$H^2$ are CR diffeomorphisms by the
assumption that their Segre homomorphisms are injective.
We should also point out that if
$M$ is of D-finite type (finite type in the sense of
D'Angelo [D'A]) at $0$, then it does not contain a nontrivial
formal variety through $0$. 

The study of formal mappings between analytic objects has a
long history. We shall mention here only results which
are related to this work. In their study of normal
forms for real-analytic Levi nondegenerate hypersurfaces in
complex space, Chern and Moser [CM] proved that any formal
invertible mapping between such is convergent. We should
also mention that Moser and Webster [MW] proved convergence
of formal equivalences between real-analytic two dimensional
surfaces with complex elliptic tangents (of a certain
kind) at the origin in $\bC^2$. On the other
hand, Moser and Webster in the same paper also proved that
there are  real-analytic surfaces through
the origin in $\bC^2$ which are formally equivalent but not
biholomorphically equivalent at $0$. Subsequently, Webster
[W] proved that each real Lagrangian surface in $\bC^2$
with a nondegenerate complex tangent at $0$ is formally
equivalent to a certain quadric. Later, Gong [Go2] showed
that there exist such real-analytic Lagrangian surfaces
which are not biholomorphic to that quadric. An analogous
situation arises in the classification of glancing
hypersurfaces (see Melrose [Mel], Oshima [O], and also Gong
[Go1]). 

Previous results on convergence of formal mappings between
real-analytic submanifolds in complex space were obtained by 
the authors in [BER3] and [BER5], where analytic
dependence of holomorphic mappings on their jets of a
prescribed finite order was studied. We should point out
that the results in the present paper are new even
for invertible mappings. (See also
the closing remark at the end of this paper.)

Finite determination (by
their 2-jets) for formal invertible mappings between smooth
Levi nondegenerate hypersurfaces was proved by Chern and Moser
[CM] (see also earlier work of E. Cartan [C] and Tanaka [Ta]).
Later work was done by e.g.\ Tumanov and Henkin [TH],
Beloshapka [B], and the authors [BER2], [BER5]. For a more
detailed history, the reader is referred to [BER5]. 

This paper is organized as follows. In \S1, we introduce the
basic concepts of formal manifolds and mappings in the setting
of ideals in the ring of formal power series. We also give
the fundamental properties of the Segre mappings for generic
submanifolds: these mappings play an important role in
the proofs. The precise formulations of our main results are
given in \S2. Another ingredient in the proofs is the
reflection identity for formal mappings presented in \S3. The
proofs of the results stated above and those in \S2 are then
given in
\S4--7.

\subhead Acknowledgements\endsubhead The authors would like to
thank J. D'Angelo, X. Gong, and X. Huang for helpful
bibliographical information, and also C. Huneke for his help
with the proof of Lemma 3.32.

\heading \S1. Preliminaries and basic
definitions\endheading

\subhead \S1.1. Formal mappings and
manifolds\endsubhead Let
$\bC\dbl x\dbr =\bC\dbl x_1,\ldots,x_k\dbr $ be the ring of
formal power series in $x=(x_1,\ldots, x_k)$
with complex coefficients. We shall say that 
a proper ideal
$I\subset\bC\dbl x_1,\ldots, x_k\dbr $ is a
{\it manifold ideal} of dimension $l$, if it is
generated by
$k-l$ power series $h_1,\ldots, h_{k-l}$ whose
differentials at the origin are independent. 
This is equivalent to the local ring
of
$I$,
$\bC\dbl x\dbr /I$, being a regular local ring of
Krull dimension $l$ (see e.g. [Ha], [AM]). We shall
write $\Cal I(h_1(x),\ldots, h_{k-l}(x))=\Cal
I(h(x))$ for the ideal in $\bC\dbl x\dbr $ generated by
$h(x)=(h_1(x),\ldots, h_{k-l}(x))$.

 By a formal  mapping
$F\:(\bC^k,0)\to(\bC^m,0)$ 
 we shall mean a vector of formal power
series $F\in\bC\dbl x_1,\ldots x_k\dbr ^m$ such that each
component of
$F=(F_1,\ldots, F_m)$ has no constant
term. Writing $\bC\dbl y\dbr = \bC\dbl y_1,\ldots,
y_m\dbr $, we let $\eta_F\:\bC\dbl y\dbr \to
\bC\dbl x\dbr $ 
be the ring homomorphism defined by 
$$
\eta_F(g)(x)=(g\circ F)(x),\quad g\in\bC\dbl y\dbr .
\tag1.1.1
$$
Conversely, any ring holomorphism $\eta: \bC\dbl y\dbr \to
\bC\dbl x\dbr $ is of the form \thetag{1.1.1} for some
uniquely determined formal mapping $F$. Given an ideal
$I\subset\bC\dbl x\dbr $, we define its {\it pushforward}
by $F$,
$F_*(I)\subset \bC\dbl y\dbr $, as
follows
$$
F_*(I)=\{f\in\bC\dbl y\dbr \:f\circ F\in I\}.\tag1.1.2
$$
Note that $F_*(I)=\eta_F^{-1}(I)$, with $\eta_F$
defined by \thetag{1.1.1}. 

It is easy to check that $F_*(I)$ is an ideal
and if $I$ is
prime, then $F_*(I)$ is also prime. However, 
even if $I$ is a manifold ideal,
$F_*(I)$ need not be a manifold ideal. Observe
that the mapping
$\eta_F$ is a ring isomorphism if and only if $m=k$
and the Jacobian of $F$ at $0$ is invertible.  In that
case we shall refer to $F$ as a {\it formal change of
coordinates}.  If $I \subset \bC\dbl x\dbr $ is an ideal
then $I$ is a manifold ideal of dimension $l$ if and
only if there is a formal change of coordinates $F$
such that
$F_*(I )\subset \bC\dbl y\dbr $ is generated by 
the coordinates $y_1,\ldots, y_{k-l}$.

A {\it formal vector field} in $\bC^k$ is a derivation
of
$\bC\dbl x\dbr $.  It is easy to check that any formal
vector field
$X$ can be written uniquely in the form
$$
X=\sum_{j=1}^ka_j(x)\frac{\partial}
{\partial x_j},\tag1.1.3
$$
with $a_j\in\bC\dbl x\dbr $.  We shall say that $X$ 
 is tangent to a proper ideal $I\subset\bC\dbl x\dbr $ if
for any
$f\in I$ we have $Xf\in I$. If $F$ is a formal change
of coordinates as above, then we may define the
{\it pushforward} of a formal vector field $X$ by
$$
(F_*X) (g) = \eta_{F^{-1}}(X(\eta_F(g))), \ \quad g
\in
\bC\dbl y\dbr .
\tag 1.1.4
$$
It is clear that $F_*X$ is also a formal vector
field.   Furthermore, if $X$ is tangent to an ideal
$I\subset \bC\dbl x\dbr $, then $F_*X$ is tangent to
$F_*(I)$.  

If
$Z=(Z_1,\ldots, Z_N)$,
$\zeta=(\zeta_1,\ldots,\zeta_N)$, and $I\subset
\bC\dbl Z,\zeta\dbr $ an ideal, we shall say that $I$ is a
{\it real ideal} if  $f(Z,\zeta) \in I$ implies
$\bar f(\zeta, Z) \in I$.  Here $\bar f$ is the
formal power series obtained by taking the complex
conjugates of the coefficients of $f$. If
$f(Z,\zeta) \in \bC\dbl Z,\zeta\dbr $, we shall say that
$f$ is {\it real} if $f(Z,\zeta) =\bar f(\zeta, Z)
$ The notion of a real ideal is not invariant
under  all formal changes of coordinates in
$\bC\dbl Z,\zeta\dbr $, and we shall restrict ourselves to
 changes of coordinates of the form $Z'=H(Z)$,
$\zeta'=\bar H(\zeta)$, where $H(Z)$ is a formal
change of coordinates in $\bC^N$.  We shall write
$\Cal H(Z, \zeta) :=(H(Z),\bar H(\zeta))$. If $I$ is a
real ideal in $\bC\dbl Z,\zeta\dbr $, so is the ideal $\Cal
H_* (I)$, as defined by \thetag{1.1.2}. 

If $I \subset
\bC\dbl Z,\zeta\dbr $ is a real manifold ideal of
dimension $2N-d$ then one can find generators
$\rho=(\rho_1,\ldots, \rho_d)
$ which are real and have
linearly independent differentials at the origin,
i.e. satisfy   
$$
\partial\rho_1(0)\wedge
\ldots\wedge \partial\rho_d(0)\neq 0.\tag1.1.5
$$
If $\rho'=(\rho'_1,\ldots,\rho'_d)$
is another set of generators of $I$ (not necessarily
real), then there exists a $d\times d$ matrix of
formal power series
$a(Z,\zeta)$ (necessarily invertible at
$0$) such that
$$
\rho(Z,\zeta)\sim
a(Z,\zeta)\rho'(Z,\zeta).\tag1.1.6
$$ 
We say that a real manifold ideal
$I\subset\bC\dbl Z,\zeta\dbr $ of dimension $2N-d$ defines
a 
 {\it formal real submanifold} $M$ of
$\bC^N$ through
$0$ of codimension $d$ (and dimension $2N-d$), and
we write $I = \Cal I(M)$. We shall say that a
formal vector field in $\bC^N\times\bC^N$ is
tangent to $M$ if it is tangent to $\Cal I(M)$.

The motivation behind these definitions is the
following. If in addition the components of the
formal series
$\rho$ above are convergent, then the
equations
$\rho(Z,\bar Z)=0$ define a real-analytic
submanifold $M$ of $\bC^N$ through $0$.
 Also, if $M$ is a smooth real
submanifold in $\bC^N$ through $0$, then the
Taylor series at $0$ of a smooth defining
function $\rho$ of $M$ near $0$ defines a formal
real submanifold of $\bC^N$ through $0$. Thus,
the notion of a formal real submanifold is a
generalization of that of real-analytic and
smooth real submanifolds of $\bC^N$. 

A formal real submanifold $M$ of $\bC^N$ is said
to be {\it generic} if there exist
formal  series $\rho$ generating $\Cal I(M)$ 
satisfying the condition 
$$
\partial_Z\rho_1(0)\wedge
\ldots\wedge\partial_Z\rho_d(0)\neq 0,\tag1.1.7
$$
which in particular implies \thetag{1.1.5}. 

If $H\:(\bC^N,0)\to (\bC^{N'},0)$ is a formal mapping,
then as before we associate to it a formal
mapping
$\Cal H\:(\bC^N\times\bC^N,0)\to
(\bC^{N'}\times\bC^{N'},0)$ defined by 
$$
\Cal H(Z,\zeta)=(H(Z),\bar H(\zeta)).\tag 1.1.8
$$
If $M$ and $M'$ are formal real submanifolds through $0$ in $\bC^N$ and 
$\bC^{N'}$, respectively, then we say that the formal mapping $H$ maps 
$M$ into $M'$, denoted $H\:(M,0)\to (M',0)$, if $\Cal
I(M')\subset
\Cal H_*(\Cal I(M))$. (When $M$ and $M'$ are
real-analytic and $H$ is  convergent, i.e. defines a
holomorphic mapping near $0$, then $H(M)
\subset M'$ if and only if the formal mapping 
$H$ maps the formal submanifold $M$ into the
formal submanifold
$M'$.) If $N=N'$ and 
$H$ is invertible, then the ring homomorphism $\eta_{\Cal H}\:
\bC\dbl Z',\zeta'\dbr \to\bC\dbl Z,\zeta\dbr $, defined by
\thetag{1.1.1}, where
$(Z',\zeta')\in \bC^{N'}\times
\bC^{N'}$ and $(Z,\zeta)\in \bC^N\times \bC^N$,
is an isomorphism. If, in addition, $\dim  
M=\dim M'$ and $H\:(M,0)\to (M',0)$, then $\Cal
H_*(\Cal I(M))=\Cal I(M')$ and $\eta_{\Cal H}$
is an isomorphism between $\Cal I(M)$ and
$\Cal I(M')$. In this case, we say that $M$ and $M'$
are {\it formally equivalent} and that $H$ is a
{\it formal equivalence} between $M$ and $M'$. 

In this paper, we shall consider formal
mappings which are more general than formal
equivalences. For a formal generic submanifold
$M$ of codimension $d$ through the origin in $\bC^N$, we
define the {\it Segre variety ideal} $I_1(M)\subset
\bC\dbl Z\dbr$ as the ideal generated by
$\rho_1(Z,0),\ldots, \rho_d(Z,0)$, where the 
$\rho_j(Z,\zeta)$ are generators of the ideal $\Cal
I(M)$\footnote{In the convergent
case, the variety defined by $I_1(M)$ in $\bC^N$ is the
so-called Segre variety of $M$ at $0$.}. The ideal $I_1(M)$ is a manifold
ideal which only depends on the ideal
$\Cal I(M)$ and not on the choice of generators
$\rho_j(Z,\zeta)$. The reader can easily check that if
$H\:(\bC^N,0)\to (\bC^{N'},0)$ is a formal mapping
which maps $M$ into another formal generic submanifold
$M'$ though the origin in
$\bC^{N'}$, then in the notation introduced
above, the homomorphism $\eta_{H}$ maps
$I_1(M')$ into $I_1(M)$. Hence, $\eta_{H}$
induces a homomorphism 
$$
\Phi_H\:\bC\dbl Z'\dbr/I_1(M') \to \bC\dbl
Z\dbr/I_1(M).\tag 1.1.9
$$
The homomorphism $\Phi_H$ is called the {\it Segre
homomorphism} of $H$. When $H$ is a formal equivalence,
this homomorphism is an isomorphism. We shall consider
in this paper formal mappings whose Segre homomorphisms
are injective. In particular, all formal finite mappings
satisfy this condition (see Lemma 6.1).

 A formal vector field in $\bC^N\times\bC^N$ of
the form
$$
X=\sum_{j=1}^N\left(a_j(Z,\zeta)\frac{\partial}
{\partial Z_j}+b_j(Z,\zeta)\frac{\partial}
{\partial \zeta_j}\right),\tag1.1.10
$$ is said to be of {\it type
$(0,1)$} if $a_j\sim0$, $j=1,\ldots, N$, and
similarly of {\it type $(1,0)$} if the $b_j\sim
0$. Let $M$ be a formal generic real
submanifold in
$\bC^N$ of codimension $d$. We let
$\Cal D_M^{0,1}$ denote the
$\bC\dbl Z,\zeta\dbr $-module generated by all formal
$(0,1)$ vector fields tangent to $M$. A formal vector
field 
$L\in \Cal D_M^{0,1}$ is usually called a {\it CR
vector field tangent to $M$}.  We define
$\Cal D^{1,0}_M$ in a similar fashion. We define the
$\bC\dbl Z,\zeta\dbr $-module
$\Cal D_M$  by $\Cal D_M=\Cal
D^{1,0}_M\oplus\Cal D_M^{0,1}$, and define $\frak
g_M$ to be the Lie algebra generated by $\Cal D_M$.
We shall say that a
collection $L_1,\ldots , L_n\in\Cal D_M^{0,1}$, with
$n=N-d$, is a {\it basis for the CR vector fields} on
$M$ if
$\pi(L_1),\ldots, \pi(L_n)$ form a basis for the
quotient module $\Cal D^{0,1}_M/\Cal I(M)\Cal
D^{0,1}_M$, where $\pi$ is the canonical projection
$\pi\: \Cal D^{0,1}_M\to \Cal D^{0,1}_M/\Cal I(M)\Cal
D^{0,1}_M$. Loosely speaking, $L_1,\ldots, L_n$
form a basis for the CR vector fields on $M$ if they
form a basis for $\Cal D^{0,1}_M$  modulo those formal
vector fields whose coefficients belong to $\Cal I(M)$.
(One can always find a basis for the CR
vector fields on $M$; see \S1.2.4 below.)

Let
$T'_0\bC^{2N}$ denote the holomorphic
tangent space of $\bC^{2N}$ at $0$, i.e.\ the
space of all tangent vectors of the form
$$
\sum_{j=1}^N\left(a_j\frac{\partial}
{\partial Z_j}\bigg |_0+b_j\frac{\partial}
{\partial \zeta_j}\bigg 
|_0\right),\tag1.1.11
$$ 
with  $a_j,b_j\in \bC$. We denote by
$\Cal D_M(0)$ and
$\frak g_M(0)$ the
subspaces of $T'_0\bC^{2N}$ obtained by
evaluating   at $0$ the coefficients of the
formal vector fields in $\Cal D_M$ and
$\frak g_M$ respectively. We say that $M$ is of {\it
finite type} at
$0$ if $\dim_\bC \frak g_M(0)=\dim M=2N-d$. 
(Note that the vector space $\Cal D_M(0)$ has
dimension $2N-2d$;  this follows easily from 
the fact that $M$ is generic and of
codimension $d$.) 

Another notion which will be used in
this paper is that of essential finiteness at $0$
for a formal generic submanifold $M$ of codimension $d$ in
$\bC^N$. Let
$L_1,\ldots, L_n$ be a basis for the CR vector fields
on $M$ as described above. We write
$$
L_j=\sum _{k=1}^N
a_{jk}(Z,\zeta)\frac{\partial}{\partial
\zeta_k},\tag 1.1.12
$$
where $a_{jk}\in\bC\dbl Z,\zeta\dbr$, and set
$$
X_j:=\sum _{k=1}^N
a_{jk}(0,\zeta)\frac{\partial}{\partial
\zeta_k},\tag 1.1.13
$$
for $j=1,\ldots, n$. Given generators
$\rho_1(Z,\zeta),\ldots, \rho_d(Z,\zeta)$ for $\Cal
I(M)$, we define, for $l=1,\ldots, d$ and
$\alpha\in\Bbb Z_+^n$,
$$
c_{l\alpha}(Z):=X^\alpha\rho_l(Z,\zeta)\big
|_{\zeta=0},\tag 1.1.14
$$
with $X^\alpha:=X_1^{\alpha_1}\ldots X_n^{\alpha_n}$.
We say that $M$ is {\it essentially finite} at $0$ if
the ideal generated by all the
$c_{l\alpha}(Z)$ has finite codimension in $\bC\dbl
Z\dbr$. The reader can check that this definition
depends only on the ideal $\Cal I(M)$ and not on the
choice of generators or basis of CR vector fields. (See
[BER4, Chapter XI] for a similar definition for a
smooth generic submanifold.) For a real-analytic generic
submanifold, the notion of essential finiteness given above
coincides with the standard one (see e.g.\ [BER4, Chapter
XI]), since an ideal $\Cal J$ in $\bC\{ Z\}$ is of finite
codimension if and only if the ideal generated by $\Cal J$ in
$\bC\dbl Z\dbr $ is of finite codimension (see e.g.\ [BER4,
Corollary 5.1.26]). Also, if $M$ is real-analytic and does not
contain a nontrivial holomorphic subvariety through $0$, then
$M$ is essentially finite at $0$ (see [BER4, Chapter XI.4]). In
fact, the analogous statement in the formal setting is also
true. The proof of this relies on Lemma 3.32 below, and will
be given in \S6 below.

\subhead \S1.2. Normal coordinates and the Segre
mappings\endsubhead We keep the notation introduced in
the previous section; in particular, $M$ is a formal generic
submanifold of codimension
$d$ with real manifold ideal $\Cal I=\Cal I(M)$. There is a
formal change of coordinates in
$\bC^N$,
$Z=Z(z,w)\in\bC\dbl z,w\dbr ^N$, with $z=(z_1,\ldots,
z_n)$, $w=(w_1,\ldots, w_d)$, and $N=n+d$, such
that if we make the corresponding change of
coordinates $\zeta=\bar
Z(\chi,\tau),$ with $ \chi=(\chi_1,\ldots,
\chi_n),\, \tau= (\tau_1,\ldots \tau_d)$, then
the ideal $\Cal I$ of
$M$, after making the above change of
coordinates, is generated by
$$
w_j-Q_j(z,\chi,\tau),\quad j=1,\ldots,d,\tag1.2.1
$$
where the
$Q_j(z,\chi,\tau)\in\bC\dbl z,\chi,\tau\dbr $ satisfy
$$
Q_j(0,\chi,\tau)\sim
Q_j(z,0,\tau)\sim\tau_j.\tag1.2.2
$$
(See [BER4, Chapter IV].) Such coordinates will be
called {\it normal coordinates}. It follows
from the reality of $\Cal I$ that
$\Cal I$ is also generated by 
$$
\tau_j-\bar Q_j(\chi,z,w),\quad
j=1,\ldots,d.\tag1.2.3
$$
In such coordinates, we can take
$$
L_j=\frac{\partial }{\partial\chi_j}+\sum_{l=1}^d \bar
Q_{l,\chi_j}(\chi,z,w)\frac{\partial}{\partial
\tau_l},\quad j=1,\ldots, n,\tag 1.2.4
$$
as a basis of the CR
vector fields on $M$. The reader can easily check that
if we write
$$
\bar Q_j(\chi,z,w)=\sum_{\alpha\in\Bbb Z_+^n}\bar
q_{j\alpha}(z,w)\chi^\alpha,\tag1.2.5
$$
then 
$M$ is essentially finite at $0$ if and only if the
ideal generated by $\{\bar
q_{j \alpha}(z,0)\}$, $j=1,\ldots d$, $\alpha\in\Bbb
Z_+^n$, has finite codimension in $\bC\dbl z\dbr$.

We shall introduce some more notation. Let $M$ and $M'$ be
formal generic submanifolds through the origin in $\bC^N$ and
$\bC^{N'}$ of codimension $d$ and $d'$, respectively. We let
$Z=(z,w)$,
$\zeta=(\chi,\tau)$, with
$z=(z_1,\ldots, z_n)$, $w=(w_1,\ldots,
w_d)$, $\chi=(\chi_1,\ldots, \chi_n)$ and
$\tau=(\tau_1,\ldots,
\tau_d)$, be formal normal coordinates for $M$, and similarly
$Z'=(z',w')$, $\zeta'=(\chi',\tau')$ be formal normal
coordinates for
$M'$. Thus, the ideal $\Cal I(M)$ is generated by
\thetag{1.2.1}, and similarly for $\Cal I(M')$
with $Q$ replaced by $Q'$. If $Z'=H(Z)$ is a formal
holomorphic mapping $(\bC^N,0)\to
(\bC^{N'},0)$, sending $M$ into $M'$, then we write
$H=(F,G)$, with $F=(F_1,\ldots, F_{n'})$ and $G=(G_1,\ldots,
G_{d'})$. The Segre
homomorphism $\Phi_H\:\bC\dbl Z'\dbr/I_1(M')\to \bC\dbl
Z\dbr/I_1(M)$ can be identified with the homomorphism
$$\tilde\Phi\:\bC\dbl z'\dbr \to \bC\dbl z\dbr,$$ where
$\tilde\Phi(h)(z)=h(F(z,0))$. Thus, injectivity of the
Segre homomorphism $\Phi_H$ is equivalent to the statement
that there is no formal nontrivial power series $h\in\bC\dbl
z'\dbr$ such that $h(F(z,0))\sim 0$. In particular, as
mentioned in
\S1.1, if $H$ is a finite formal mapping (with $N=N'$ and
$d=d'$), then $\Phi_H$ is injective 
(see Lemma 6.1). 
 
We now introduce another important ingredient in this
paper. Let $Q(z,\chi,\tau)$ be the $\bC^d$-valued power series
whose components are the $Q_j(z,\chi,\tau)$.  Consider, for
each integer
$k\geq1$, the formal mapping $v^k\:(\bC^{kn},0)\to(\bC^{N},0)$
defined as follows. For $k=2j$, 
$$
\multline
v^{2j}(z,\chi^1,\ldots,z^{j-1},\chi^j):=\bigg(z,
Q\big(z,\chi^1,\bar
Q\big(\chi^1,z^{1},Q\big(z^{1},
\chi^{2},\ldots,\\\bar
Q\big(\chi^{j-1},z^{j-1},Q\big(z^{j-1},
\chi^j,0\big)\big)\ldots\big)
\big)\big)\bigg),
\endmultline\tag1.2.6
$$  
and, for $k=2j+1$, 
$$
\multline
v^{2j+1}(z,\chi^1,\ldots,z^{j-1},\chi^j,
z^j):=\bigg(z,Q\big(z,\chi^1,\bar
Q\big(\chi^1,z^{1},Q\big(z^{1},\chi^{2},
\ldots,\\ Q\big(z^{j-1},\chi^j,\bar
Q\big(\chi^j,z^j,
0\big)\big)\ldots\big)\big)\big)\bigg),
\endmultline\tag1.2.7
$$
where $j=1,2,\ldots$. For $k=0,1$, we define
$$
v^0:=(0,0),\quad v^1(z):=(z,Q(z,0,0))=(z,0).
$$
We shall refer to the mapping $v^k$ as the
$k$th {\it Segre mapping} of $M$. 

\proclaim{Proposition 1.2.8} Let $f\in\Cal
I(M)$. Then, for any $k\geq 0$,
$$
f(v^{k+1}(z,\chi^1, z^1,\ldots),\bar
v^k(\chi^1,z^1,\ldots))\sim 0.\tag1.2.9
$$
\endproclaim

\demo{Proof} For simplicity, we only
consider the case where $k=2j$. It follows
from
\thetag{1.2.6} and
\thetag{1.2.7} that
$$
v^{2j+1}(z,\chi^1,\ldots,z^{j-1},\chi^j,
z^j)\sim\left (z,Q(z,\bar
v^{2j}(\chi^1,\ldots,z^{j-1},\chi^j,
z^j)\right).\tag1.2.10
$$
It suffices to show \thetag{1.2.9} for the
generators of $\Cal I(M)$ given by
\thetag{1.2.1}, for which \thetag{1.2.9}
is an immediate consequence of
\thetag{1.2.10}.\qed 
\enddemo

In this paper, we shall not need the explicit form of the
Segre mappings $$v^k\:(\bC^{kn},0)\to (\bC^N,0)$$ given in
\thetag{1.2.6} and \thetag{1.2.7}, but merely the properties
given by Proposition 1.2.8 and Theorem 1.2.11 below (and the
fact that the $v^k$ are convergent when $M$ is real-analytic).
It will be convenient to write
$v^{k}(z,\xi)$,
$z\in\bC^n$ and
$\xi\in \bC^{kn}$. Thus, for
fixed $k$, we write $\xi=(\chi^1,z^1,\ldots)\in\bC^{(k-1) n}$.
The equation \thetag{1.2.9} can then be written as
$$
f(v^{k+1}(z,\xi),\bar v^k(\xi))\sim 0,\tag 1.2.9'
$$ 
where, in \thetag{1.2.9'}, we have $\xi\in\bC^{kn}$.

For a formal power series mapping
$F\:(\bC^\ell,0)\to (\bC^m,0)$ we denote by
$\rk\,(F)$ the rank of the
matrix
$(\partial F_j/\partial x_k)$, $j=1,\ldots, m$, $k=1,\ldots,
\ell$, in
$\Bbb K^\ell$; here, $\Bbb K$ denotes
the quotient field of $\bC\dbl x\dbr $.  We shall need the
following result.

\proclaim{Theorem 1.2.11} Let $M$ be a
formal generic submanifold of $\bC^N$
through $0$ of codimension $d$.  Then, $M$ is of finite type
at
$0$ if and only if there exists $k_1\leq d+1$ such that
$\rk(v^k)=N$ for
$k\geq k_1$. Moreover, if $M$ is real-analytic and of finite
type at $0$, then there exists
$(z_0,\xi_0)\in\bC^n\times\bC^{(2k_1-1)n}$ arbitrarily close
to the origin such that $v^{2k_1}(z_0,\xi_0)=0$ and the rank of
the Jacobian matrix $\partial v^{2k_1}/\partial (z,\xi)$ at
$(z_0,\xi_0)$ is $N$.
\endproclaim

The proof of the first part of Theorem 1.2.11 can be found in
[BER5] (see also [BER1] and [BER4] for the case where
$M$ is real-analytic). The proof of the last statement follows
from [BER4, Proposition 10.6.19].

\head \S 2. Statement of further results; Examples\endhead 

We state now the main results of this paper from which the
theorems in the introduction will be deduced. We begin with a
result which will be shown to imply Theorem 3.

\proclaim{Theorem 2.1} Let $M$ and $M'$ be real-analytic
generic submanifolds through the origin in $\bC^N$ and
$\bC^{N'}$, respectively. Assume that
$M$ is of finite type at
$0$ and that $M'$ is essentially finite at $0$. If
$H\:(\bC^N,0)\to (\bC^{N'},0)$ is a formal holomorphic mapping
sending $M$ into $M'$ whose Segre homomorphism $\Phi_H$ is
injective, then the formal series $H$ converges in a
neighborhood of the origin. 
\endproclaim 

We shall now give a result for hypersurfaces which in
particular will imply Theorem 1.

\proclaim{Theorem 2.2} Let $M$ and $M'$ be real-analytic
hypersurfaces through the origin in $\bC^N$, $N\geq 2$. Assume
that
$M$ is essentially finite at $0$ and that $M'$ does not contain
any nontrivial holomorphic subvariety through $0$. Then any
formal mapping
$H\:(\bC^N,0)\to (\bC^N,0)$ sending $M$ into $M'$ is
convergent.   
\endproclaim 

\remark{Remark $2.3$} The condition that $M'$ above does not
contain a nontrivial holomorphic subvariety is necessary for
the conclusion of Theorem 2.2 to hold. Indeed, if $M'$ is any
real-analytic submanifold containing a nontrivial holomorphic
subvariety through the origin in
$\bC^N$, then there exists a formal mapping $H\:(\bC^N,0)\to
(\bC^N,0)$ sending all of $\bC^N$ into $M'$ (even into the
holomorphic subvariety contained in $M'$) which does not
converge. To see this, let $s\mapsto A(s)$ be a
nontrivial holomorphic mapping from a neighborhood of the
origin in the complex plane into the holomorphic subvariety
contained in $M'$ with $A(0)=0$. The reader can easily verify
that there exists $f\in\bC\dbl s\dbr$ with $f(0)=0$ such that
at least one component of $A\circ f$ is not convergent. We may
then take $H(Z)=A(f(Z_1))$. \endremark\medskip

As a byproduct of the proof of Theorem 2.2, we obtain
the following reflection principle for CR mappings
between real-analytic hypersurfaces. 

\proclaim{Theorem 2.4} Let $M$ and $M'$ be real-analytic
hypersurfaces through the origin in $\bC^N$, $N\geq 2$. Assume
that
$M$ is essentially finite at $0$ and $M'$ does not contain any
nontrivial holomorphic subvariety through $0$. Then every
smooth CR mapping $h\:(M,0)\to (M',0)$ which does not vanish
to infinite order at $0$ necessarily extends as a holomorphic
mapping $(\bC^N,0)\to (\bC^N,0)$ sending $M$ into $M'$. 
\endproclaim

In the case $N=2$, this result follows from the more general
theorem of Huang [Hu] in which the only condition assumed on
the CR mapping is continuity. However, for $N\geq 3$, Theorem
2.4 is new, since previous results in this context assumed more
stringent conditions on the CR mapping. The reflection
principle for real-analytic submanifolds in complex space has
a long history; the reader is referred to the notes in [BER4,
Chapter IX] for further information. 

We return to the case of higher codimensional generic
submanifolds. Theorem 5 will be a consequence of
the following result on finite determination, or more
precisely its corollary given below.

\proclaim{Theorem 2.5} Let $M$ and $M'$ be formal 
generic submanifolds through the origin in $\bC^N$ and
$\bC^{N'}$, respectively. Assume that $M$ is of finite type at
$0$, that $M'$ is essentially finite at $0$, and that there
exists a formal holomorphic mapping $H\:(\bC^N,0)\to
(\bC^{N'},0)$ sending $M$ into $M'$ whose Segre homomorphism
$\Phi_{H}$ is injective. Then there exists an
integer $K$ with the following property. If
$H'\:(\bC^N,0)\to (\bC^{N'},0)$ is a formal
holomorphic mapping sending $M$ into $M'$ and if 
$$
\frac{\partial^{|\alpha|} H'}{\partial Z^\alpha}(0)= 
\frac{\partial^{|\alpha|} H}{\partial
Z^\alpha}(0),\quad\forall |\alpha|\leq K,\tag 2.6
$$
then $H\sim H'$. 
\endproclaim

\proclaim{Corollary 2.7} Let $M$  be a formal 
generic submanifold through the origin in $\bC^N$. Assume 
that $M$ is of finite type and essentially finite at
$0$. Then there exists an integer $K$ with the following
property. If $M'$ is a formal 
generic submanifold through the origin in $\bC^N$ of the same
dimension as $M$, and if $H^1, H^2\:(\bC^N,0)\to
(\bC^N,0)$ are formal invertible
mappings sending $M$ into $M'$ which satisfy 
$$
\frac{\partial^{|\alpha|} H^1}{\partial Z^\alpha}(0)= 
\frac{\partial^{|\alpha|} H^2}{\partial
Z^\alpha}(0),\quad\forall |\alpha|\leq K,\tag 2.8
$$
then $H^1\sim H^2$.
\endproclaim

\demo{Proof of Corollary $2.7$} We claim that it suffices to
take
$K$ to be the integer given by Theorem 2.5 with $M'=M$ and
$H(Z):=Z$. To see this, let $M'$, $H^1$, $H^2$ be as in
Corollary 2.7 and observe that if \thetag{2.8} holds then
$\partial^\alpha ((H^1)^{-1}\circ H^1)(0)=\partial
^\alpha((H^1)^{-1}\circ H^2)(0)$, for all $|\alpha|\leq K$. By
Theorem 2.5 and the choice of $K$, we deduce that
$((H^1)^{-1}\circ H^2)(Z)\sim Z$ and, hence, the conclusion of
Corollary 2.7.\qed
\enddemo

We conclude this section by giving three simple examples
motivating each of the three conditions in Theorem 2.1.

\remark{Example $2.9$} Let $M=M'\subset \bC^3$ be the
generic submanifold of codimension two defined by
$$
\im w_1=|z|^2,\quad \im w_2=0.\tag 2.10
$$
Observe that $M$ is essentially finite but not of finite type
at $0$. Let $f\in\bR\dbl w_2\dbr$ be nonconvergent with
$f(0)=f'(0)=0$. Then, the formal mapping $H\:(\bC^3,0)\to
(\bC^3,0)$ defined by
$$
H(z,w_1,w_2):=(z,w_1,w_2+f(w_2))\tag 2.11
$$
is invertible, sends $M$ into $M$, but does not converge.
\endremark\medskip

\remark{Example $2.12$} Let $M=M'\subset \bC^4$ be the
generic submanifold of codimension two defined by
$$
\im w_1=|z_1z_2|^2,\quad \im w_2=|z_1z_2|^4.\tag 2.13
$$
In this case, $M$ is of finite type but not essentially finite
at $0$. Let $f\in\bC\dbl z_1\dbr$ be
nonconvergent with
$f(0)=0$. Then, the formal mapping $H\:(\bC^4,0)\to
(\bC^4,0)$ defined by
$$
H(z_1,z_2,w_1,w_2):=(z_1e^{f(z_1)},z_2e^{-f(z_1)},w_1,w_2)\tag
2.14
$$
is invertible, sends $M$ into $M$, but does not converge.
\endremark\medskip

\remark{Example $2.15$} Let $M=M'\subset \bC^4$ be the
generic submanifold of codimension two defined by
$$
\im w_1=|z_1|^2-|z_2|^2,\quad \im w_2=|z_1|^4-|z_2|^4.\tag
2.16 
$$
This time, $M$ is essentially finite and of finite type 
at $0$. Let $f\in\bC\dbl z_1\dbr$ be
nonconvergent with
$f(0)=0$. Then, the formal mapping $H\:(\bC^4,0)\to
(\bC^4,0)$ defined by
$$
H(z_1,z_2,w_1,w_2):=(f(z_1),f(z_1),0,0)\tag
2.17
$$
sends $M$ into $M$, but does not converge. Observe that the
Segre homomorphism
$\Phi_H$ is not injective.
\endremark\medskip

\head \S 3. A reflection identity\endhead

The following reflection identity will play an important role
in the proofs of Theorems 2.1 and 2.5.

\proclaim{Theorem 3.1} Let $M$ and $M'$ be formal 
generic submanifolds through the origin in $\bC^N$  and
$\bC^{N'}$, respectively. Assume that $M'$ is essentially
finite at
$0$. Let
$H:(\bC^N,0)\to (\bC^{N'},0)$ be a formal holomorphic mapping
sending $M$ into
$M'$ whose Segre homomorphism $\Phi_H$ is injective. Then
there exist a positive integer
$r$, positive integers $N_j$,
$j=1,\ldots, N'$, and monic polynomials in $X$ of the form 
$$
P_j(X,Z,\zeta,(a_\gamma)_{|\gamma|\leq
r}):= X^{N_j}+\sum_{k=0}^{N_j-1}
c^{jk}(Z,\zeta,(a_\gamma)_{|\gamma|
\leq r},)X^k, \tag 3.2
$$ 
for $j=1,\ldots, N'$ such that 
$$
P_j(H_j(Z),Z,\zeta,(\partial^\gamma \bar
H(\zeta)-\partial^\gamma\bar H(0))_{|\gamma|\leq r})\in
\Cal I(M).\tag 3.3
$$
Here,
$\gamma\in\Bbb Z^N_+$,
$a_\gamma=(a_{\gamma m})_{1\leq m\leq {N'}}$, and the $c^{jk}$
are formal power series whose coefficients depend only
on $M$, $M'$ and on the values
$\partial^\beta 
H(0)$ for $|\beta|\leq r$. In addition, if $H'$
is another formal mapping $(\bC^N,0)\to (\bC^{N'},0)$ sending
$M$ into $M'$
and if $\partial^\beta H'(0)=\partial ^\beta H(0)$,
for
$|\beta|\leq r$, then \thetag{3.3} holds with $H$
replaced by $H'$.

Moreover, if $M$ and $M'$ are real-analytic, then each of the
formal power series in \thetag{3.2} defines a holomorphic
function in a neighborhood of the origin.
\endproclaim

In order to prove Theorem 3.1, 
we let $Z=(z,w)$,
$\zeta=(\chi,\tau)$, with
$z=(z_1,\ldots, z_n)$, $w=(w_1,\ldots,
w_d)$, $\chi=(\chi_1,\ldots, \chi_n)$ and
$\tau=(\tau_1,\ldots,
\tau_d)$, be formal normal coordinates for $M$, and similarly
$Z'=(z',w')$, $\zeta'=(\chi',\tau')$ be formal normal
coordinates for
$M'$, as in \S1.2. We shall first prove Theorem 3.1
in these coordinates.
We let $L_1,\ldots, L_n$ be the basis for the formal CR vector
fields on $M$ given by \thetag{1.2.4}. Let
$H=(F,G)\:(\bC^{n}\times\bC^{d},0)\to
(\bC^{n'}\times\bC^{d'},0)$ be a formal holomorphic mapping
sending $M$ into $M'$. The following proposition gives the
desired polynomial identities for the components of $F$. 

\proclaim{Proposition 3.4} Let $M$, $M'$, and $H$ be as
above. Assume that $M'$ is essentially finite at $0$, and
that the Segre homomorphism $\Phi_H$ is injective. Then there
exist a positive integer $r$, positive integers $N_j$,
$j=1,\ldots, n'$, and monic polynomials 
$$
P_j(X,(a_\gamma)_{|\gamma|\leq r},(b_\gamma)_{|\gamma|\leq
r}):= X^{N_j}+\sum_{k=0}^{N_j-1} c^{jk}((a_\gamma)_{|\gamma|
\leq r},(b_\gamma)_{|\gamma|\leq r})X^k, \tag 3.5
$$ 
for $j=1,\ldots, n'$ such that 
$$
P_j(F_j(z,w),(L^\gamma \bar
F(\chi,\tau)-L^\gamma\bar F(0))_{|\gamma|\leq r},(L^\gamma
\bar G(\chi,\tau))_{|\gamma|\leq r})\in
\Cal I(M).\tag 3.6
$$
Here,
$\gamma\in\Bbb Z^n_+$,
$a_\gamma=(a_{\gamma m})_{1\leq m\leq {n'}}$,
$b_\gamma=(b_{\gamma l})_{1\leq l\leq {d'}}$, and the $c^{jk}$
are formal power series whose coefficients depend only
on $M'$ and on the values
$(L^\beta 
\bar F)(0)$ for $|\beta|\leq r$. In addition, if $H'=(F',G')$
is another formal mapping $(\bC^N,0)\to (\bC^{N'},0)$ sending
$M$ into $M'$
and if $L^\beta \bar F'(0)=L^\beta \bar F(0)$, for
$|\beta|\leq r$, then \thetag{3.6} holds with $(F,G)$
replaced by $(F',G')$.

Moreover, if $M$ and $M'$ are real-analytic, then each of the
formal power series in \thetag{3.5} defines a holomorphic
function in a neighborhood of the origin.
\endproclaim

For the proof of Proposition 3.4, we shall need the
following preliminary results. 

\proclaim{Lemma 3.7} Let $p_1(X,Y,Z)$ and $p_2(X,Y,Z)$
be formal power series of the form
$$
\align
p_1(X,Y,Z) &=X^N+\sum_{j=0}^{N-1} a_j(Y,Z) X^j,\tag 3.8\\
p_2(X,Y,Z) &=Y^M+K(X,Y,Z),\tag 3.9
\endalign
$$
where $X,Y\in\bC$, $Z=(Z',Z'')\in\bC^{k'}\times
\bC^{k''}$, $a_j(0,0)=0$, and $K(X,Y,Z',0)\sim 0$. Then the
ideal $\Cal I(p_1,p_2)\subset
\bC\dbl X,Y,Z\dbr$ contains a power series of the form
$$
r(Y,Z)=Y^{MN^2}+K'(Y,Z),\tag
3.10
$$
with $K'(Y,Z',0)\sim 0$. 
\endproclaim 

\demo{Proof} For $\xi=(\xi_1,\ldots,\xi_N)$, we denote
by $\sigma_0(\xi),\ldots,\sigma_{N-1}(\xi)$ the
usual elementary symmetric polynomials of $\xi$ defined by 
$$\multline
\sigma_{N-1}
(\xi) =-\sum_{j=1}^N\xi_j,\qquad
\sigma_{N-2}(\xi) =\sum_{1\leq j<k\leq N}\xi_j\xi_k,\qquad
\ldots,\\ 
\sigma_0(\xi) =(-1)^N\prod_{j=1}^N\xi_j.\endmultline\tag 3.11
$$
Consider the ring
$$
\Cal R=\bC\dbl X,Y,Z,\xi\dbr/\Cal
J,\tag
3.12
$$
where $\Cal J$ denotes
the ideal generated by $\sigma_j(\xi)-a_j(Y,Z)$, 
$j=0,\ldots, N-1$, where the $a_j$ are the coefficients of
$p_1$ in \thetag{3.8}. We use the notation
$[f]\in\Cal R$ for the equivalence class of a power series
$f\in\bC\dbl X,Y,Z,\xi\dbr$, and we also consider $\bC\dbl
X,Y,Z\dbr$ as a subring of $\bC\dbl
X,Y,Z,\xi\dbr$. Observe that
$$
[p_1]=\bigg [\prod_{j=1}^N(X-\xi_j)\bigg].\tag 3.13
$$
Consider the power series
$$
\hat q(Y,Z,\xi)=\prod_{j=1}^{N}p_2(\xi_j,Y,Z).\tag 3.14
$$
Since $\hat q$ is symmetric in the $\xi$, it follows from
Newton's theorem for symmetric functions (see e.g.\ [VW])
that there exists a power series
$q(Y,Z,\sigma_0,\ldots,\sigma_{N-1})$ such that
$$
\hat q(Y,Z,\xi)=q(Y,Z,\sigma_0(\xi),\ldots
\sigma_{N-1}(\xi)).\tag 3.15
$$
By the definition of the ideal $\Cal J$, we have 
$$
[\hat q(Y,Z,\xi)]=[q(Y,Z,a_0(Y,Z),\ldots
a_{N-1}(Y,Z))].\tag 3.16
$$
We claim that $[q(Y,Z,a_0(Y,Z),\ldots
a_{N-1}(Y,Z))^N]$ is in the ideal generated by $[p_1]$ and
$[p_2]$ in $\Cal R$. This is straightforward to check in view
of
\thetag{3.13--16} since $\hat q^N$ is in the ideal generated
by
$p_2$ and $\prod_{j=1}^N(X-\xi_j)$ in $\bC\dbl
X,Y,Z,\xi\dbr$. (The latter fact follows easily from the
observation that, for any $j$,  $\hat q$ is in the ideal
generated by
$p_2$ and $X-\xi_j$.)

It follows from the claim above that there are power series
$c_j(X,Y,Z,\xi)$, $j=1,2$, such that
$$
\multline
q(Y,Z,a_0(Y,Z),\ldots
a_{N-1}(Y,Z))^N-c_1(X,Y,Z,\xi)p_1(X,Y,Z)-\\
c_2(X,Y,Z,\xi)p_2(X,Y,Z) \in\Cal J.
\endmultline\tag 3.17
$$
We put $u(Y,Z):=q(Y,Z,a_0(Y,Z),\ldots
a_{N-1}(Y,Z))$. It is easy to see that $u(Y,Z)$ is of the form
$Y^{MN}+K''(Y,Z)$ with $K''(Y,Z',0)\sim 0$.
We shall complete the proof of Lemma 3.7 by
showing that $r(Y,Z):=u(Y,Z)^N$ is in the ideal $\Cal
I(p_1,p_2)$. Since the generators of the ideal
$\Cal J$ are invariant under permutations of the components of
$\xi$, we may assume without loss of generality that $c_j$,
$j=1,2$, in
\thetag{3.17} are symmetric in the $\xi$. Thus, again by
Newton's theorem, there are power series
$d_j(X,Y,Z,\sigma)$, $j=1,2$,
$\sigma=(\sigma_0,\ldots,\sigma_{N-1})$, such that   
$$
\multline
r(Y,Z)-d_1(X,Y,Z,\sigma_0(\xi),\ldots,\sigma_
{N-1}(\xi))p_1(X,Y,Z)-\\ d_2(X,Y,Z,\sigma_0(\xi),\ldots,
\sigma_ {N-1}(\xi))p_2(X,Y,Z)
\in\Cal J.
\endmultline\tag 3.18
$$
It follows that
$$
\multline
r(Y,Z)-d_1(X,Y,Z,a_0(Y,Z),\ldots,a_
{N-1}(Y,Z))p_1(X,Y,Z)-\\ d_2(X,Y,Z,a_0(Y,Z),\ldots,
a_ {N-1}(Y,Z))p_2(X,Y,Z)
\in\Cal J.
\endmultline\tag 3.19
$$ 
Since the power series in \thetag{3.19} is independent of
$\xi$, the desired fact that $r(Y,Z)$ is in $\Cal I(p_1,p_2)$
will follow from the injectivity of the canonical 
homomorphism $$\bC\dbl X,Y,Z\dbr\to\Cal R.$$ 
Thus, we must
show that if $g(X,Y,Z)\in \Cal J$, then $g\sim 0$. By
expanding in
$X$ and using the fact that the generators of $\Cal J$ are
independent of $X$, it suffices to prove that if
$$
f(Y,Z)\sim\sum_{j=0}^{N-1}e_j(Y,Z,\xi)
(\sigma_j(\xi)-a_j(Y,Z)),\tag 3.20
$$
for some power series $e_j$, then $f\sim 0$. By considering
all possible formal curves $t\mapsto (Y(t),Z(t))$ through
the origin, we are reduced to proving the following
statement: If 
$$
h(t)\sim\sum_{j=0}^{N-1}\tilde e_j(t,\xi)
(\sigma_j(\xi)-a_j(Y(t),Z(t))),\tag 3.21
$$
for some power series $\tilde e_j$, then $h\sim 0$. By the
formal Puiseux expansion (see e.g. [BK]), it
follows that there exists an integer $J$ and formal series
$f_1(w),\ldots, f_N(w)$ such that
$$
p_1(X,Y(w^J),Z(w^J))\sim\prod_{j=1}^N(X-f_j(w)).\tag 3.22
$$
Making the substitutions $t=w^J$ and $\xi_j=f_j(w)$,
$j=1,\ldots, N$, in \thetag{3.21}, we conclude that
$h(w^J)\sim 0$ and hence $h\sim 0$. This completes the proof
of Lemma 3.7. \qed\enddemo 

We shall make use of Lemma 3.7 to prove the following
proposition.

\proclaim{Proposition 3.23} Let
$f(u,v)=(f_1(u,v),\ldots,f_r(u,v))$, with $u\in\bC^p$ and
$v\in\bC^q$, be formal power series with vanishing constant
terms. Assume that the ideal generated by $f_j(0,v)$,
$j=1,\ldots, r$, has finite codimension in $\bC\dbl v\dbr$.
Then, there exist power series $P_j(u,v)$, $j=1,\ldots, q$, of
the form
$$
P_j(u,v)=v_j^{N_j}+\sum_{k=0}^{N_j-1} b_{jk}(u)v_j^k,\tag
3.24
$$
where the $b_{jk}$ are power series in $u$ with vanishing
constant terms, such that each $P_j(u,v)$ is in 
the ideal
$\Cal I(f(u,v))\subset \bC\dbl u,v\dbr$.
\endproclaim
\demo{Proof} We
first reduce the situation to a simpler one. By standard
arguments (see e.g.\  [BER4, Proposition 5.1.5]), there are
power series
$a_{jk}(v)$ and an integer
$N$ such that
$$ v_j^N\sim\sum_{k=1}^r
a_{jk}(v)f_k(0,v),\quad j=1,\ldots,q.\tag 3.25
$$ 
Define power series 
$g_{kl}(u,v)$ such that
$$
f_k(u,v)\sim f_k(0,v)+\sum_{l=1}^pu_l\, g_{kl}(u,v),\quad
k=1,\ldots, r.\tag 3.26
$$
It follows that
$$ 
v_j^N+K_j(u,v)\sim \sum_{k=1}^r
a_{jk}(v)f_k(u,v),\tag 3.27
$$ 
where
$$ 
K_j(u,v):=\sum_{k=1}^r
\sum_{l=1}^p a_{jk}(v)u_l\, g_{kl}(u,v),\quad
j=1,\ldots,q.\tag 3.28
$$
Observe that $v_j^N+K_j(u,v)$ are in the ideal $\Cal
I(f(u,v))$ and that $K_j(0,v)\sim 0$. 

Hence to prove Proposition
3.23 we may assume that $q=r$ and the $f_j$ are of the form 
$$
f_j(u,v):=v_j^N+K_j(u,v),\quad j=1,\ldots q,\tag 3.29
$$
where $K_j(0,v)\sim 0$. We reason by induction on $q$. 
For
$q=1$, the desired conclusion follows from the (formal)
Weierstrass Preparation Theorem (see e.g.\ [ZS]) applied to the
power series 
$f_1(u,v)$, where $v=v_1$. We
shall now show how to reduce the case of $q$
to that of $q-1$.

We apply the Weierstrass Preparation Theorem to
$f_1(u,v)$. Hence there exist
power series 
$c_l(u,v_2,\ldots,v_q)$, $l =
0,\ldots,N-1$ satisfying
$c_l(0,v_2,\ldots,v_n) \sim 0$ so that
$$  
f_1(u,v)\sim U(u,v)\left(v_1^{N} +
\sum_{l=0}^{N-1}c_l(u,v_2,\ldots,v_q)v_1^l\right),
\tag 3.30
$$ 
where $U(u,v)$ is a unit in $\bC\dbl u,v\dbr$. We apply, for
$j=2,\ldots, q$, Lemma 3.7 with $X=v_1$, $Y=v_j$, $Z''=u$,
$Z'$ being the components $v_k$ for $k\geq 2$ and $k\neq j$,  
$p_1(X,Y,Z)$ being $v_1^{N} +
\sum_{l=0}^{N-1}c_l(u,v_2,\ldots, v_q)v_1^l$, and
$p_2(X,Y,Z)$ being $f_j(u,v)$. We
conclude that the ideal $\Cal I(f(u,v))$
contains power series of the form
$$
f'_j(u,v_2,\ldots, v_q):= v_j^{N^3}+K'_j(u,v_2,\ldots,
v_{q}),\quad j=2,\ldots q,\tag 3.31
$$
with $K'_j(0,v_2,\ldots,v_q)\sim 0$. By the inductive
hypothesis, we conclude that $\Cal I(f(u,v))$
contains power series $P_j(u,v_j)$, $j=2,\ldots, q$, of the
form \thetag{3.24}. We obtain $P_1(u,v_1)$ by repeated
application of Lemma 3.7 in a similar fashion as above.
This completes the proof of Proposition 3.23.\qed\enddemo 

For the proof of Proposition 3.4, we also need the following
result which can be viewed as a formal version of the
Nullstellensatz.

\proclaim{Lemma 3.32} Given $K(x)=(K_1(x),\ldots, K_r(x))$
with 
$K_j(x)\in\bC\dbl x\dbr$, 
$x\in\bC^n$, the following are equivalent:
\roster
\item"(i)" $\dim_\Bbb C \bC\dbl x\dbr/\Cal I(K(x))=
\infty$,
\item"(ii)" there exist $\mu(s)=(\mu_1(s),...,\mu_n(s))$,
with 
$\mu_1,...,\mu_n\in
\bC\dbl s\dbr$ and $s\in\bC$, such that
$\mu(0)=0$, $\mu(s)\not\sim 0$, and
$$
K_j(\mu(s))\sim0,\qquad\forall\, 1\leq j\leq r.\tag 3.33
$$
\endroster\endproclaim
\remark {Remark} The authors are indebted to C. Huneke 
for the proof of the implication
(i)$\Rightarrow$(ii).\endremark\medskip 

\demo{Proof} {\bf (i)$\Rightarrow$(ii)}. Since the ideal 
$\Cal I(K(x))$ has  infinite (vector space) codimension
(i.e. (i) in the lemma is satisfied),  there is a prime ideal
$\p$ containing
$\Cal I(K(x))$ that also has infinite codimension (see
e.g.\ [BR, Lemma 3.3]). Since $\frak m$, the maximal ideal in
$\bC\dbl x\dbr$, is the only proper
prime ideal of finite codimension, we have
$\p\subsetneqq\m$. We take $\p$ to be maximal with this
property,  i.e. $\Cal I(K(x))\subset
\p\subsetneqq\m$ and there is no prime ideal 
$\q$ (of infinite codimension) such that
$\p\subsetneqq\q\subsetneqq\m$. Then the  Krull
dimension of $\hO/\p$ is precisely 1 ([ZS, p. 218]). The
integral closure of $\hO/\p$  in its field of
fractions is then a local 1-dimensional integrally closed
ring, i.e.\ a complete discrete valuation ring denoted by $R$
below, whose residue field,  being finite over
$\Bbb C$, is precisely $\Bbb C$. It then follows from  the
Cohen Structure Theorem (see e.g. [ZS, p. 307]) that R is
isomorphic to $\bC\dbl s\dbr$, with $s\in\bC$. Thus, we obtain
a homomorphism 
$\tilde\psi\:\hO/\Cal I(K(x))
\rightarrow\bC\dbl s\dbr$
by composing the homomorphism $\hO/\Cal I(K(x))\rightarrow 
\hO/\p$, induced by the inclusion $\Cal
I(K(x))\hookrightarrow\p$, with the inclusion
$\hO/\p\hookrightarrow R\cong \bC\dbl s\dbr$. Let
$\psi\:\hO\rightarrow\bC\dbl s\dbr$ be the homomorphism
obtained by composing the projection $\hO\rightarrow\hO/\Cal
I(K(x))$ with
$\tilde\psi$. By construction, this homomorphism is not
identically 0 and $\Cal I(K(x))\subset\ker \psi$. The desired
formal power series map $\mu(s)$ is given by
$\mu_j(s)=\psi(x_j)$, which are not all 0 since $\psi$ is not
identically 0. To verify
\thetag{3.33} it suffices to note  that, for any
$h\in\hO$, we have
$\psi(h)(s)=h(\mu(s))$. The identity
\thetag{3.33} follows from the fact that $\Cal
I(K(x))\subset\ker\psi$.
\smallskip
\noindent{\bf (ii)$\Rightarrow$(i)}. Assume, in order  to
reach a contradiction, that $\Cal I(K(x))$ has finite
codimension. Then there are integers $N_j$, for $1\leq j\leq
n$, such that $x_j^{N_j}$ is in $\Cal I(K(x))$, i.e. there are
$A^k_j\in\hO$ ($1\leq j\leq n$ and $1\leq k\leq r$) such that
$$
x_j^{N_j}\sim\sum_{k=1}^r A^k_j(x)K_k(x).\tag 3.34
$$
Substituting $x=\mu(s)$ into this we obtain, using 
\thetag{3.33},
$$
\mu_j(s)^{N_j}\sim\sum_{k=1}^r
A_j^k(\mu(s))K_k(\mu(s))\sim0,\tag 3.35
$$
for each $1\leq j\leq n$.
This contradicts the fact that the $\mu_j(s)$ are not all 0, 
proving the desired implication.\qed\enddemo

\demo{Proof of Proposition $3.4$}  
The fact that $H$ sends $M$ into $M'$ is equivalent to 
$$
\bar G(\chi,\tau)-\bar Q'(\Bar
F(\chi,\tau),F(z,w),G(z,w))\in\Cal I(M)^{d'}.\tag 3.36
$$
In what follows, we write $F=F(z,w)$, $\bar F=\bar
F(\chi,\tau)$ and similarly for $G$ and $\bar G$. We write 
$$
\bar G - \bar Q'(\bar F,F,G) = \bar G-
\bar Q'(\bar F,F,0)-P(\bar F,F,G)G\in\Cal I(M)^{d'}, \tag 3.37
$$ 
where $P$ is a $d'\times d'$ matrix of  
formal power series in $2n'+d'$ variables.  We decompose
$$ 
\bar Q'(\bar F,F,0) = \sum_\alpha
\bar Q'_\alpha(F)\bar F^\alpha. 
\tag 3.38
$$  

We claim that the ideal $\Cal J \subset \bC\dbl \mu\dbr$ 
generated by
$\sum_{\alpha}\bar Q'_\alpha(\mu)L^\beta\bar F^\alpha(0)$,
$\beta\in\Bbb Z_+^n$, has
finite codimension. To see this, first observe that 
$$
L^\beta h(0):=L^\beta
h(\chi,\tau)|_0=\frac{\partial^{|\beta|}
h}{\partial\chi^\beta}(0), \tag 3.39
$$
for any power series $h\in\bC\dbl\chi,\tau\dbr$. Now, pick an
integer $r$ such that the ideal $\Cal J$ is generated by the
power series $\sum_{\alpha}\bar Q'_\alpha(\mu)L^\beta\bar
F^\alpha(0)$ for
$|\beta|\leq r$. Recall that, in normal coordinates,
injectivity of the Segre homomorphism $\Phi_H$ is equivalent
to the fact that there are no nontrivial power series
$h\in\bC\dbl \chi'\dbr$ such that $h(\bar F(\chi,0))\sim 0$.
Hence, there are
no nontrivial power series $g\in\bC\dbl s,\chi'\dbr$ such that
$g(s,\bar F(\chi,0))\sim 0$.  It follows that for any formal
curve
$s\mapsto
\mu(s)$ with $\mu=(\mu_1,\ldots, \mu_n)$, $\mu(0)=0$, and
$\mu\not\sim 0$, we have 
$$
\sum_{\alpha}\bar Q'_\alpha(\mu(s))\bar
F(\chi,0)^\alpha\not\sim 0,\tag 3.40
$$
because otherwise $\bar Q'_\alpha(\mu(s))\sim 0$ for all
$\alpha\in\Bbb Z_+^{n'}$ which would contradict the essential
finiteness of
$M'$ by the implication (ii)$\implies$(i) in Lemma 3.32. (We
have used here the characterization of essential finiteness in
normal coordinates given in \S1.2.) It follows by using
\thetag{3.39} that, given a formal curve
$\mu(s)$ as above, there exists $\beta\in\Bbb Z_+^n$ such that
$$
\sum_{\alpha}\bar Q'_\alpha(\mu(s))L^\beta\bar
F^\alpha(0)\not\sim 0.\tag 3.41
$$
By the choice of $r$ above, we may assume that $|\beta|\leq
r$. The claim follows by applying the implication
(i)$\implies$ (ii) in Lemma 3.32 to the collection of power
series 
$\sum_{\alpha}\bar Q'_\alpha(\mu)L^\beta\bar F^\alpha(0)$,
$|\beta|\leq r$.

We now apply
$L^{\beta}$, $|\beta|\leq r$, to \thetag{3.37} and
substitute $G=Q'(F,\bar F,\bar G)$ to obtain
$$ 
\sum_\alpha
\bar
Q'_\alpha(F)L^{\beta}\bar F^\alpha+L^
{\beta}\big(P(\bar F,F,G)\big)G\big|_{G=Q'(F,\bar F,\bar
G)}-L^{\beta}\bar G\in
\Cal I(M)^{d'}.
\tag 3.42
$$ 
Observe that $L^\beta\bar G(0)=0$, by normality of the
coordinates, but $L^\beta\bar F(0)$ is in general not 0. We
can rewrite \thetag{3.42} as
$$ 
\multline
\sum_\alpha
\bar
Q'_\alpha(F)(L^{\beta}\bar F^\alpha-L^\beta\bar
F(0)^\alpha)+\sum_\alpha
\bar
Q'_\alpha(F)L^{\beta}\bar F(0)^\alpha+\\ 
S\left(F,\bar F,\bar G, (L^\gamma\bar F-L^\gamma\bar
F(0))_{1\leq|\gamma|\leq r}\right)Q'(F,\bar F,\bar
G)-L^{\beta}\bar G\in
\Cal I(M)^{d'},\endmultline
\tag 3.43
$$ 
for some formal power series $S$ which depends on
$M'$ and the values $L^\gamma\bar F(0)$ for $|\gamma|\leq
r$. (We have used here the fact that $L^
{\beta}\big(P(\bar F,F,G)\big)$ is polynomial in the variables
$L^\gamma\bar F$ for
$1\leq |\gamma|\leq r$.) Let $a_\gamma=(a_{\gamma p})_{1\leq
p\leq {n'}} = (L^\gamma\bar F_p-L^\gamma\bar F_p(0))_{1\leq
p\leq {n'}}$ and
$b_\gamma=(b_{\gamma l})_{1\leq l\leq {d'}} = (L^\gamma\bar
G_l)_{1\leq l\leq {d'}}$, considered as independent
variables, for $|\gamma|\leq r$.  We write $a =
(a_{\gamma p})$ and $b =
(b_{\gamma l})$.
The equations \thetag{3.43} can be
written in the form
$$ 
R_{\beta j}((L^\gamma\bar F_p-L^\gamma\bar
F_p(0)),(L^\gamma\bar G_l),F)\in \Cal I(M),\ \
j=1,\ldots,{d'},\ |\beta|\leq r,\tag 3.44
$$ 
where the
$R_{\beta j}$ are formal power series whose coefficients
depend only on $M'$ and the values $L^\beta\bar F(0)$ for
$|\beta|\leq r$. Furthermore, using the normality of the
coordinates, we observe that $R_{\beta}(0,0,F):=(R_{\beta
j}(0,0,F))_{1\leq j\leq d'}\sim
\sum_\alpha
\bar
Q'_\alpha(F)L^{\beta}\bar F(0)^\alpha$. Thus, the ideal
generated by $R_{\beta j}(0,0,\mu)$, $j=1,\ldots d'$ and
$|\beta|\leq r$, has finite codimension by the claim above.
The first part of Proposition 3.4 then follows from
Proposition 3.23 applied to the power series $R_{\beta
j}(a,b,F)$ (where $(a,b,F)$ are considered as independent
variables). Moreover, if $H'=(F',G')$ is another
formal mapping sending
$M$ into $M'$ and if $L^\beta\bar F'(0)= L^\beta\bar F(0)$,
for $|\beta|\leq r$, then the ideal generated by the
power series $\sum_\alpha\bar Q'(\mu)L^\beta\bar F'(0)$, for
$|\beta|\leq r$, has finite codimension in $\bC\dbl\mu\dbr$.
The fact that $H'$ satisfies the same polynomial identity as
$H$ follows from the construction above. 

To prove the last statement of Proposition 3.4, i.e. when
$M$ and $M'$ are assumed to be real-analytic, we observe that
both Proposition 3.23 and Lemma 3.32 (which are used in
the proof above) hold in the convergent setting. Indeed, for
the proof of Proposition 3.23 in the convergent case the
reader is referred to [BER4, Theorem 5.3.9]. Lemma 3.32 in
the convergent case is a special case of the Nullstellensatz
for convergent power series (see e.g.\ [Gu, Theorem
II.E.2]). The convergence of the power series in \thetag{3.5}
now follows by an inspection of the proof above in the formal
case. This completes the proof of Proposition 3.4.
\qed\enddemo

\demo{Proof of Theorem $3.1$} To complete the proof of Theorem
3.1, using Proposition 3.4, we observe that
$$
G(z,w)-Q'(F(z,w),\bar F(\chi,\tau,)\bar G(\chi,\tau))\in\Cal
I(M)^{d'}.\tag 3.41
$$
We obtain polynomial identities for the components $G$,
similar to those for $F$, by repeated applications of
Lemma 3.7 as in the proof of Proposition 3.23. To obtain
the identities in \thetag{3.3}, it suffices to replace the
vector fields $L_1,\ldots, L_n$ in the identities
\thetag{3.6} for
$F$, and the similar ones just obtained for 
$G$, by their expressions in the coordinates
$(z,w,\chi,\tau)$. The fact that Theorem 3.1 holds in any
set of coordinates follows from the result in normal
coordinates by an application of Proposition 3.23. This
completes the proof of Theorem 3.1.\qed\enddemo

\heading \S. 4. Proof of Theorem 2.1\endheading

For the proof of Theorem 2.1, we shall make use of Theorem
3.1  and the following consequence of the approximation theorem
of Artin  (see [A]).
\proclaim{Theorem 4.1} Let $h\in\bC\dbl x\dbr$,
$x=(x_1,\ldots, x_m)$, assume that $P(h(x),x)\sim 0$, where 
$$
P(X,x)=\sum_{j=0}^J a_j(x)X^j,\quad a_j\in\bC\{x\},\
a_J(x)\not\equiv 0,\tag 4.2
$$
then $h\in\bC\{x\}$. Here, $\bC\{x\}$ denotes the ring of
convergent power series in $x$.\endproclaim

\demo{Proof of Theorem $4.1$} By Artin's approximation theorem,
for every integer $k$, there exists a convergent solution
$X=\tilde h_k(x)$ of $P(X,x)=0$ such that the power series of
$\tilde h_k(x)$ agrees with that of $h(x)$ up to order $k$.
The conclusion in Theorem 4.1 is a consequence of Lemma 4.3
below whose proof is left to the reader.\qed\enddemo

\proclaim{Lemma 4.3} Let $P(X,x)$ be of the form
$$
P(X,x)=\sum_{j=0}^J a_j(x)X^j,\quad a_j\in\bC\dbl x\dbr,\
a_J(x)\not\equiv 0.\tag 4.4
$$
Then there exists a positive integer $m(P)$ such that if 
$h,h'\in\bC\dbl x\dbr$ satisfy 
$$
\aligned
&P(h(x),x)\sim P(h'(x),x)\sim
0\\ &\partial ^\alpha h(0)=\partial^\alpha h'(0),
\quad \forall |\alpha|\leq m(P),\\ 
\endaligned \tag 4.5
$$
then $h(x)\sim h'(x)$. 
\endproclaim

\demo{Proof of Theorem 2.1}
Consider the following vector fields tangent to $M$ near the
origin
$$
S_j=\frac{\partial}{\partial Z_j}-\rho_{Z_j}(Z,\zeta)
(\rho_{\zeta''}(Z,\zeta))^{-1}\frac{\partial}{\partial
\zeta''},\quad j=1,\ldots N,\tag 4.6
$$
where $\rho=(\rho_1,\ldots,\rho_d)$ is a defining
$\bR^d$-valued function for $M$ near $0$ with 
$\rho_{\zeta''}=(\partial\rho_j/\partial \zeta_{k})_{1\leq
j,k\leq d}$ invertible at $0$ (which we may assume without
loss of generality); in \thetag{4.6}, $\partial/\partial
\zeta''=(\partial/\partial \zeta_k)_{1\leq k\leq d}$ is
considered as a $d\times 1$ matrix and
$\rho_{Z_j}=(\rho_{k,Z_j})_{1\leq k\leq d}$ as a $1\times
d$ matrix. Observe that $S^\alpha h=\partial ^{|\alpha|}
h/\partial Z^\alpha $ for any $h\in \bC \dbl Z\dbr$. 

We shall make use of the Segre mappings of $M$ introduced in
\S1.2. Recall the notation $v^{k+1}(z,\xi)$, $z\in\bC^n$
and
$\xi\in\bC^{(k-1)n}$, for the $(k+1)$th Segre mapping of $M$ at
$0$. (We use the notation $v^0$ for the constant
mapping $v^0=0$.) Consider the following property, for
$k\geq 0$ and
$\alpha\in\Bbb Z_+^N$,
$$
\aligned 
&(\partial^\beta H)\circ
v^l\ \text{{\rm is convergent}},\quad  \forall 0\leq l\leq
k-1,\
\forall \beta\in\Bbb Z_+^N,\\ &(\partial^\beta H)\circ
v^k\ \text{{\rm is convergent}},\quad  \forall\beta\leq \alpha,
\endaligned\leqno(*)_{k,\alpha}
$$
where $\beta<\alpha$ means that $\beta$ precedes $\alpha$ in
the lexicographical ordering of $\Bbb Z_+^N$. We shall say that
$(*)_k$ holds if
$(*)_{k,\alpha}$ holds for all $\alpha\in\Bbb Z_+^N$. We shall
prove
$(*)_k$ for all
$k$ by induction. Observe that $(*)_0$ holds, since
$v^0$ is the constant mapping. We assume $(*)_{k}$, and
wish to prove $(*)_{k+1}$. We first prove $(*)_{k+1,0}$. We
shall make use of Theorem 3.1. We substitute
$(Z,\zeta)=(v^{k+1}(z,\xi),\bar v^k(\xi))$ in \thetag{3.3}. By
Proposition 1.2.8, we have for
$j=1,\ldots, N'$
$$
P_j(H_j(
v^{k+1}(z,\xi)),v^{k+1}(z,\xi),\bar v^k(\xi),(\partial^\gamma
\bar H(\bar v^k(\xi))-\partial^\gamma\bar H(0))_{|\gamma|\leq
r})\sim 0.\tag 4.7
$$
It follows from the inductive hypothesis $(*)_k$ that 
$$
P_j(X,v^{k+1}(z,\xi),\bar v^k(\xi),(\partial^\gamma
\bar H(\bar v^k(\xi))-\partial^\gamma\bar H(0))_{|\gamma|\leq
r})\in\bC\{z,\xi\}[X].\tag 4.8
$$ 
Hence, by Theorem 4.1 and \thetag{4.7}, we conclude that
$(*)_{k+1,0}$ holds. To complete the proof of $(*)_{k+1}$, we
prove
$(*)_{k+1,\alpha}$ for all $\alpha$ by induction on $\alpha$
(using the lexicographic ordering of $\Bbb Z_+^N$).  For this
we shall need the following result whose proof will be given
later. 
Recall that if $\Cal J\subset \bC\dbl x\dbr$ is an ideal, then
we say that a formal vector field $S$ is tangent to $\Cal J$
if $S$ is a derivation of $\bC\dbl x\dbr$ preserving $\Cal
J$. 

\proclaim{Lemma 4.9} Let $\Cal J\subset \bC\dbl x\dbr$,
$x=(x_1,\ldots, x_m)$, be an ideal and $S_1,\ldots, S_N$
formal vector fields tangent to $\Cal J$. Assume that
$h\in\bC\dbl x\dbr$, and that there exists a formal power
series $P(X,x)$ of the form
$$
P(X,x)=\sum_{j=0}^J a_j(x)X^j,\quad a_j\in\bC\dbl
x\dbr,\ a_J(x)\sim 1,\tag 4.10
$$
with $P(h(x),x)\in\Cal J$. Then, for any formal mapping
$v\:(\bC^r,0)\to (\bC^m,0)$ such that $f\circ v\sim 0$ for
all $f\in\Cal J$ and any $\alpha\in\Bbb Z_+^N$, there exists a
formal power series $R_\alpha(X,x)$ of the form
$$
R_\alpha(X,x)=\sum_{j=0}^J b_{\alpha j}(x)X^j,\quad
b_{\alpha j}\in\bC\dbl x\dbr,\tag 4.11
$$
with $R_\alpha(X,v(y))\not\sim 0$ and 
$$
R_\alpha\big((S^\alpha h)(v(y)),v(y)\big)\sim 0.\tag 4.12
$$
Here, we have used the notation $S^\alpha=S_1^{\alpha_1}\ldots
S_N^{\alpha_N}$. Furthermore, the coefficients $b_{\alpha j}$
are universal polynomials in $S^\beta a_j$, $|\beta|\leq
|\alpha|J$ and
$j=0,\ldots, J-1$, and $S^\gamma h$, for $\gamma<\alpha$,
where $\gamma<\alpha$ means that $\gamma$ precedes $\alpha$
in the lexicographic ordering of
$\Bbb Z_+^N$.  
\endproclaim

We proceed with the proof of Theorem 2.1.
Assume
$(*)_{k+1,\alpha^0}$ for some $\alpha^0\in\Bbb Z_+^N$. Recall
that
$(*)_{k+1,0}$ has already been proved. We complete the
induction by proving $(*)_{k+1,\alpha'}$, where $\alpha'$ is
the multi-index immediately following $\alpha^0$ in the
lexicographic ordering of $\Bbb Z_+^N$.

By applying Lemma
4.9 to the polynomials
$$
P_j(X,Z,\zeta,(\partial^\gamma \bar
H(\zeta)-\partial^\gamma\bar H(0))_{|\gamma|\leq r}), \quad
j=1,\ldots, N',
$$
with $\Cal J=\Cal I(M)$, $S_j$ given by
\thetag{4.6}, $x=(Z,\zeta)$,
$h=H_j$, and
$v=(v^{k+1},\bar v^k)$, we obtain, for each $j=1,\ldots, N'$
and every
$\alpha\in\Bbb Z_+^N$, power series
$R_{j\alpha}(X,Z,\zeta)$ in $\bC\dbl Z,\zeta \dbr [X]$ such
that
$R_{j\alpha}(X,v^{k+1}(z,\xi),\bar v^k(\xi))\not\sim 0$ and
$$
R_{j\alpha}(\partial^\alpha H_j(v^{k+1}(z,\xi)),
v^{k+1}(z,\xi),\bar v^k(\xi))\sim 0.\tag 4.13
$$
Moreover, the coefficients of $R_{j\alpha}$, as a polynomial
in $X$, are polynomials in $\partial ^\beta H_j(Z)$,
for
$\beta<\alpha$, and
$S^\gamma c^{jl}(Z,\zeta,(\partial^\delta\bar
H_j(\zeta)-\partial^\delta\bar H_j(0))_{|\delta|\leq r})$,
$\gamma\leq N_j\alpha$, where $c^{jl}$ are given by
\thetag{3.2}. Thus, by the
inductive hypothesis
$(*)_{k+1,\alpha^0}$, it follows that $$R_{j\alpha'}(
X,v^{k+1}(z,\xi),\bar v^k(\xi))\in\bC\{z,\xi\}[X].$$ The
property $(*)_{k+1,\alpha'}$ follows from Theorem 4.1. This
completes the proof of property $(*)_k$ for all $k$. 

We shall apply $(*)_{k,0}$ for $k=2k_1$, where $k_1$ is given
by Theorem 1.2.11. By Theorem 1.2.11, there exists
$(z^0,\xi^0)\in\bC^n\times\bC^{(2k_1-1)n}$, arbitrarily close
to $(0,0)$, such that the rank of $\partial v^{2k_1}/\partial
(z,\xi)$ at $(z^0,\xi^0)$ is $N$ and $v^{2k_1}(z^0,\xi^0)=0$.
Since $(z,\xi)\mapsto H(v^{2k_1}(z,\xi))$ is holomorphic in a
neighborhood of $(0,0)$ by $(*)_{2k_1}$, we may choose
$(z^0,\xi^0)$ in that neighborhood. By applying the implicit
function theorem, we may find a right inverse of
$v^{2k_1}$,  
$Z\mapsto
\theta(Z)$ with $\theta(0)=(z^0,\xi^0)$. Hence,
$H(v^{2k_1}(\theta(Z)))$ is a convergent mapping in a
neighborhood of $0$. This completes the proof of Theorem
2.1 modulo the proof of Lemma 4.9, since
$v^{2k_1}(\theta(Z))\equiv Z$ near
$0\in \bC^N$.

\demo{Proof of Lemma $4.9$} Observe that \thetag{4.7} holds for
$\alpha=0$ with $R_0(X,x)=P(X,x)$. We will need the following
Leibnitz formula
$$
S^\gamma
P(h(x),x)=\sum_{k=0}^J\,\sum_{\mu+\nu^1+\ldots+\nu^k=\gamma}
\frac{\gamma!}{\mu!\nu^1!\ldots\nu^k!}\, (S^\mu a_k)
(S^{\nu^1}h)\ldots (S^{\nu^k}h).\tag 4.14
$$
The reader can check that \thetag{4.14} can be rewritten as
follows. For any given
$\alpha\leq\gamma$, 
$$
\multline
S^\gamma
P(h(x),x)=\\ A (\gamma,\alpha,0)+\sum_{j=1}^J\sum_
{{\nu^1\geq\alpha,\ldots,\nu^j\geq\alpha}\atop{\gamma-\nu^1-
\ldots-\nu^j\in\Bbb Z_+^N}} 
A(\gamma,\alpha,j,\nu^1,\ldots,\nu^j) (S^{\nu^1}h)\ldots
(S^{\nu^j}h),\endmultline\tag 4.15
$$
where 
$$
\multline
A(\gamma,\alpha,j,\nu^1,\ldots,\nu^j)=\\ \sum_{k=j}^J\binom
{k}{j}\sum_{{\mu,\nu^{j+1},\ldots,\nu^k}\atop{{\mu+\nu^1+
\ldots+
\nu^{k}=\gamma}\atop
{\nu^{j+1}<\alpha,\ldots,\nu^k<\alpha}}}\frac{\gamma!}{\mu!
\nu^1!\ldots\nu^k!}(S^\mu a_k) (S^{\nu^{j+1}}h)\ldots
(S^{\nu^k}h),\endmultline\tag 4.16
$$  
$$
A(\gamma,\alpha,0)=
\sum_{k=0}^J\sum_{{\mu+\nu^1+\ldots+\nu^{k}=\gamma}\atop
{\nu^{1}<\alpha,\ldots,\nu^k<\alpha}}\frac{\gamma!}{\mu!
\nu^1!\ldots\nu^k!}(S^\mu a_k) (S^{\nu^{j+1}}h)\ldots
(S^{\nu^k}h).\tag 4.17
$$
Observe that the first term in the sum in \thetag{4.16} (i.e.
the term corresponding to $k=j$) is 
$$
\frac{\gamma!}{(\gamma-\nu^1-\ldots-\nu^j)!\nu^1!\ldots\nu^j!}
S^{\gamma-\nu^1-\ldots-\nu^j}a_j.
$$
For
$\gamma=\alpha$, we have, from \thetag{4.15},
$$
S^\alpha
P(h(x),x)=A(\alpha,\alpha,0) +
A(\alpha,\alpha,1,\alpha) S^\alpha h.\tag 4.18
$$
If $A(\alpha,\alpha,1,\alpha)\circ v\not\sim 0$, then we take
$R_\alpha(X,x)=A(\alpha,\alpha,0)+
A(\alpha,\alpha,1,\alpha) X$. If not, then we denote by
$\gamma^0$ the first multi-index $>\alpha$ for which
$A(\gamma^0,\alpha,j,\nu^1,\ldots,\nu^j)\circ v\not \sim 0$
for some $1\leq j\leq J$, $\nu^1\geq\alpha, \ldots,
\nu^j\geq \alpha$. Such a $\gamma^0$ exists (and is
$\leq J\alpha$), since
$A(J\alpha,\alpha,J,\alpha,\ldots, \alpha)\sim
(J\alpha)!/(\alpha!)^J$ as is easily verified. We claim that
necessarily
$A(\gamma^0,\alpha,j_0,\alpha,\ldots,\alpha)\circ v\not\sim 0$
for some $1\leq j_0\leq J$ and
$A(\gamma^0,\alpha,j,\nu^1,\ldots,\nu^j)\circ v\sim 0$ for
all $1\leq j\leq J$ and $\nu^1\geq\alpha,\ldots,
\nu^j\geq\alpha$ such that $\sum_{1\leq l\leq
j}\nu^l>j \alpha$. Indeed, the following identity is a
consequence of \thetag{4.16}
$$
A(\gamma^0,\alpha,j,\nu^1,\ldots,\nu^j)=e_{\gamma^0,\alpha,
j,\nu^1,\ldots,\nu^j}
A(\gamma^1,\alpha,j,\alpha,\ldots,\alpha),\tag 4.19
$$
where 
$$
\gamma^1=\gamma^0+j\alpha-\sum_{l=1}^j\nu^l,\tag 4.20
$$
and $e_{\gamma^0,\alpha,
j,\nu^1,\ldots,\nu^j}$ is a positive number depending only
on $\gamma^0,\alpha,
j,\nu^1,\ldots,\nu^j$. Observe that
$\gamma^1\geq \alpha$, and $\gamma^1<\gamma^0$ if 
$\sum_{1\leq l\leq j}\nu^l>j \alpha$. It follows from the
definition of $\gamma^0$ and \thetag{4.19} that
$A(\gamma^0,\alpha,j,\nu^1,\ldots,\nu^j)\circ v\sim 0$ for
all $1\leq j\leq J$ and $\nu^1\geq\alpha,\ldots,
\nu^j\geq\alpha$ such that $\sum_{1\leq l\leq
j}\nu^l>j \alpha$. The claim is proved. The conclusion of the
lemma follows by taking
$$
R_\alpha(X,x):= A (\gamma^0,\alpha,0)+\sum_{j=1}^J
A(\gamma^0,\alpha,j,\alpha,\ldots,\alpha)
X^j,\tag 4.21
$$
since $R_\alpha(S^{\alpha}h(x),x)|_{x=v(y)}\sim
S^{\gamma^0}  P(h(x),x)|_{x=v(y)}$ and $S^{\gamma^0} 
P(h(x),x)\in \Cal J$.
\qed\enddemo

This completes the proof of Theorem 2.1.\qed\enddemo

\heading \S 5. Proof of Theorem 2.5\endheading

The proof of Theorem 2.5 is parallel to that of Theorem
2.1, and we shall keep the notation from \S 4. 
We again use the notation $v^k(z,\xi)$, $z\in\bC^n$ and
$\xi\in\bC^{(k-1)n}$, for the $k$th (formal) Segre mapping of
$M$ at
$0$. Let $H$ be the formal mapping given in the statement of
Theorem 2.5. In what follows, $\Cal F(M,M')$ will denote the
set of all formal mappings $(\bC^N,0)\to (\bC^{N'},0)$ that
send
$M$ into $M'$. Consider the following property, for
$k\geq 0$ and
$\alpha\in\Bbb Z_+^N$.
\medskip
\flushpar
{\it $(**)_{k,\alpha}$ There exists $K(k,\alpha)\in\Bbb Z_+$
such that for any $H'\in\Cal F(M,M')$ with 
$$
\partial^\beta H(0)=\partial^\beta H'(0),\ \forall
|\beta|\leq K(k,\alpha),
\tag 5.4
$$
the following holds} 
$$
(\partial^\alpha H)\circ
v^k\sim (\partial^\alpha H')\circ
v^k.\tag 5.5
$$
\medskip

Observe that $(**)_{0,\alpha}$ holds with
$K(0,\alpha)=|\alpha|$, since $v^0$ is the constant mapping
into $0$. We shall prove that $(**)_{k,\alpha}$ holds for all
$k$ and $\alpha$ by induction on $k$ and $\alpha$. First assume
$(**)_{l,\beta}$ holds for $0\leq l\leq k$ and all
$\beta\in\Bbb Z^N_+$. We shall prove
$(**)_{k+1,0}$; i.e. we must show the existence of the
integer
$K(k+1,0)$ in $(**)_{k+1,0}$. Let
$P_j(X,Z,\zeta,(a_\gamma)_{|\gamma|\leq r})$, $j=1,\ldots, N'$,
be the polynomials in $X$ given by Theorem 3.1. Pick an
integer $\tilde K\geq 0$ and consider the formal mappings
$H'\in\Cal F(M,M')$ satisfying 
$$
\partial^\beta H(0)=\partial^\beta H'(0),\ \forall
|\beta|\leq \tilde K.
\tag 5.6
$$
By {\it a priori} requiring
$\tilde K\geq r$, where
$r$ is the integer given by Theorem 3.1, we may assume that
any
$H'\in\Cal F(M,M')$ that satisfies \thetag{5.6} also
satisfies
\thetag{3.3} with $H$ replaced by $H'$. If we also take
$\tilde K\geq K(k,\beta)$, for all $|\beta|\leq r$, then
$(**)_{k,\beta}$ implies that \thetag{5.5} holds for
$|\beta|\leq r$. It follows from the
above and Proposition 1.2.8 that
$X=H_j\circ v^{k+1}$ and
$X=H_j'\circ v^{k+1}$ are both solutions of the equation
$$
P^{k+1}_j(X,z,\xi):=P_j(X,v^{k+1}(z,\xi),\bar
v^k(\xi),(\partial^\gamma
\bar H(\bar v^k(\xi))-\partial^\gamma\bar
H(0))_{|\gamma|\leq r})\sim 0.\tag 5.7
$$
Hence, if $m(P^{k+1}_j)$ is the integer given by Lemma 4.3
and we choose
$K(k+1,0)=\max (\tilde K,m(P^{k+1}_j)_{1\leq j\leq N'})$, then
the identity $H\circ v^{k+1}\sim H'\circ v^{k+1}$ follows from
Lemma 4.3. The property
$(**)_{k+1,0}$ is proved.

We now fix an integer $k$ and a multi-index $\alpha^0$. We
complete the induction by assuming
$(**)_{l,\beta}$ for all pairs $l,\beta$ satisfying either
$0\leq l<k$ and
$\beta\in\Bbb Z_+^N$ or $l=k$ and $\beta<\alpha^0$ in the
lexicographic ordering of
$\Bbb Z_+^N$. We shall prove $(**)_{k,\alpha'}$, where
$\alpha'$ is the multi-index immediately following $\alpha^0$
in the ordering of $\Bbb Z^N_+$. Again, consider those
$H'\in\Cal F(M,M')$ that satisfy \thetag{5.6} with $\tilde
K\geq r$. As above, the components $H'_j$ of such an
$H'$ satisfy the same identities
\thetag{3.3} as $H_j$, $j=1,\ldots,N'$. By Lemma 4.9, there
exist formal series $R_{j\alpha'}(X,Z,\zeta)$, $j=1,\ldots,
N'$, in
$\bC\dbl Z,\zeta
\dbr [X]$ such that
$R_{j\alpha'}(X,v^{k}(z,\xi),\bar v^{k-1}(\xi))\not\sim 0$ and
$$
R_{j\alpha'}(\partial^{\alpha'} H_j(v^{k}(z,\xi)),
v^{k}(z,\xi),\bar v^{k-1}(\xi))\sim 0.\tag 5.8
$$
Moreover, the coefficients of $R_{j\alpha'}$, as a polynomial
in $X$, are polynomials in $\partial ^\beta H_j(Z)$,
for
$\beta<\alpha'$, and
$S^\gamma c^{jl}(Z,\zeta,(\partial^\delta\bar
H_j(\zeta)-\partial^\delta\bar H_j(0))_{|\delta|\leq r})$,
$\gamma\leq N_j\alpha'$, where $c^{jl}$ are given by
\thetag{3.2}. If $H'\in\Cal F(M,M')$ satisfies
\thetag{5.6}, then there is a corresponding formal series
$R'_{j\alpha'}(X,Z,\zeta)$ for $H'$ which is obtained from
the same equation \thetag{5.7} with $H$ replaced by $H'$.
Hence, by Lemma 4.9, if
$$
\aligned
&(\partial^\beta H)\circ
v^k\sim (\partial^\beta H')\circ
v^k,\ \beta<\alpha'\\&(\partial^\gamma H)\circ
v^{k-1}\sim (\partial^\gamma H')\circ
v^{k-1},\ |\gamma|<r+ |\alpha'|\max\, (N_j)_{j=1,\ldots N'},
\endaligned
\tag 5.9
$$
then 
both $X=\partial^{\alpha'}H\circ v^k$ and
$X=\partial^{\alpha'}H'\circ v^k$ satisfy the same equation
$$
R^k_{j\alpha'}(X,z,\xi):=R_{j\alpha'}(X,
v^{k}(z,\xi),\bar v^{k-1}(\xi))\sim 0.\tag 5.10 
$$
By the inductive hypothesis
$(**)_{l,\beta}$ for pairs
$l,\beta$ as described above,
\thetag{5.9} holds if
$\tilde K\geq K(k,\beta)$, $\beta<\alpha'$, and $\tilde K\geq
K(k-1,\gamma)$, $|\gamma|<r+ |\alpha'|\max\,
(N_j)_{j=1,\ldots N'}$. The property $(**)_{k,\alpha'}$
follows from Lemma 4.3 if we choose $K(k,\alpha')$ to be the
maximum of
$\tilde K$, as described above, and $m(R^k_{j\alpha'})$ for
$j=1,\ldots, N'$. This completes the induction, and proves
that $(**)_{k,\alpha}$ holds for all $k$ and $\alpha$. 

To complete the proof of Theorem 2.5, recall that $M$ is
of finite type at the origin and let
$k_1$ be the integer obtained by applying Theorem 1.2.11 to
$M$. If $H'\in \Cal F(M,M')$
satisfies \thetag{5.6} with $\tilde K=K(k_1,0)$ then, by
property $(**)_{k_1,0}$, we have $$H\circ v^{k_1}\sim H'\circ
v^{k_1}.\tag 5.11$$ Since $\rk(v^{k_1})=N$ the components of
the formal mapping $v^{k_1}$ cannot satisfy a nontrivial
formal relation (see e.g.\ [BER4, Proposition 5.3.5]). Hence,
by \thetag{5.11} we must have $H_j-H'_j\sim 0$. This completes
the proof of Theorem 2.5.
\qed

\heading \S6. Proofs of Theorems 3, 4, 5\endheading

For the proofs of Theorems 3, 4, and 5, we shall need the
following two lemmas.

\proclaim{Lemma 6.1} Let $M$ and $M'$ be formal generic
submanifolds of the same
dimension through the origin in $\bC^N$. Assume that
$H\:(\bC^N,0)\to (\bC^N,0)$ is a finite formal mapping sending
$M$ into $M'$. Then the Segre homomorphism $\Phi_H$ is
injective.
\endproclaim

\demo {Proof} Let $Z=(z,w)$ and $Z'=(z',w')$ be normal
coordinates for $M$ and $M'$, respectively, and write
$H=(F,G)$ as in \S1.2. Recall that injectivity of the Segre
homomorphism
$\Phi_H$ is equivalent to the fact that there is no nontrivial
formal power series $h\in\bC\dbl z'\dbr$ such that
$h(F(z,0))\sim 0$. We first claim that the formal mapping
$z\mapsto F(z,0)$ is again finite. To see this, consider the
homomorphism $\psi\:\bC\dbl z,w\dbr\to\bC\dbl z\dbr$ defined
by $\psi(f)(z)=f(z,0)$. Observe that $\psi$ is surjective and
maps the ideal $\Cal I(F(z,w),G(z,w))$ into $\Cal I(F(z,0))$
since $G(z,0)\sim 0$ by normality of the coordinates. Hence,
$\psi$ induces a surjective homomorphism 
$$ 
\tilde\psi\:\bC\dbl z,w\dbr/\Cal I(F(z,w),G(z,w))\to 
\bC\dbl z\dbr/\Cal I(F(z,0)).\tag 6.2
$$
Since the vector space on the left hand side of \thetag{6.2}
is finite dimensional and $\tilde \psi$ is surjective, we
conclude that the vector space on the right hand side of
\thetag{6.2} is also finite dimensional. This proves the
claim above. To complete the proof of the lemma, we must show
that any finite formal mapping $K\:(\bC^k,0)\to (\bC^k,0)$
induces an injective ring homomorphism $\eta_K\:\bC\dbl
x\dbr\to\bC\dbl x\dbr$. This can be proved by showing that the
Jacobian determinant $\det(\partial K/\partial x)$  cannot
vanish identically (see e.g.\ [BER4, Theorem 5.1.37]) and that
the latter implies that there is no nontrivial $h\in\bC\dbl
x\dbr$ such that $h(K(x))\sim 0$ (see [BER4, Proposition
5.3.5]).\qed 
\enddemo 

\proclaim{Lemma 6.3} Let $M$ be a formal generic submanifold
through the origin in $\bC^N$. If $M$ does not contain a
non-trivial formal variety through $0$, then $M$ is
essentially  finite at $0$. Moreover, if $M$ is real-analytic
and does not contain a nontrivial holomorphic variety through
$0$, then $M$ is essentially finite at $0$. 
\endproclaim

\demo{Proof} Let $d$ be the codimension of $M$. If $n:=N-d=0$,
then there is nothing to prove since every such generic
submanifold is essentially finite at $0$. (In this case, there
are no CR vector fields tangent to $M$ and the condition for
essential finiteness is void.) Assume $n\geq 1$ and that $M$
does not contain any nontrivial formal curve through $0$. 
Let
$Z=(z,w)$ be normal coordinates for $M$ as in \S1.2. Thus, $M$
is given by \thetag{1.2.1}. To check that $M$ is essentially
finite at $0$ we use the expansion \thetag{1.2.5} and must
show that the ideal generated by $\bar q_{j\alpha}(z,0)$,
$j=1,\ldots, d$ and 
$\alpha\in\Bbb Z_+^n$, in $\bC\dbl z\dbr$ has finite
codimension. If not, then by Lemma 3.32 there would exist a
nontrivial formal curve $\mu\:(\bC,0)\to (\bC^n,0)$ such that
$$
\bar q_{j\alpha}(\mu(s),0)\sim 0,\quad j=1,\ldots, d,\
\alpha\in\Bbb Z_+^n.\tag 6.3
$$
This would contradict the assumption that $M$ does not contain
a nontrivial formal curve through $0$, since, as the
reader can easily verify, it would follow that
$s\mapsto (\mu(s),0)$ is a nontrivial formal curve through
$0$ contained in $M$. This proves the first part of
the lemma. The proof of the last statement in
Lemma 6.3, which is more standard, is completely analogous and
uses the Nullstellensatz instead of Lemma 3.32. (See e.g.\
[BER4, Chapter XI].)\qed\enddemo 

We now give the proofs of the results in the introduction.

\demo{Proof of Theorems $3$ and $5$} By Lemmas 6.1 and 6.3,
Theorems 3 is an immediate consequence of Theorem
2.1. Similarly, Theorem 5 follows by applying Theorem 2.5 to
the formal submanifolds and the formal mapping associated to
$M$, $M'$ and $H$, respectively.\qed\enddemo

\demo{Proof of Theorem $4$} Assume that $H\:(\bC^N,0)\to
(\bC^N,0)$ is a formal invertible mapping sending $M$ into
$M'$. Since the Segre homomorphism of $H$ is injective, Theorem
4 is a consequence of Theorem 2.1 provided that we can show
that $M'$ is essentially finite at $0$. Indeed, since $M$ is
essentially finite (by Lemma 6.3) and the notion of essential
finiteness is invariant under formal invertible mappings, it
follows that $M'$ is also
essentially finite at $0$. 
\qed\enddemo

\heading \S7. Proofs of Theorems 1, 2.2, and 2.4\endheading

 For the proof of Theorem
2.2, we shall use Theorem 2.1 and the following result.

\proclaim{Proposition 7.1} Let $M$ and $M'$ be formal
hypersurfaces through the origin in $\bC^N$. Assume
that $M$ is essentially finite at $0$ and $M'$ does not
contain any nontrivial formal variety through $0$. If
$H\:(\bC^N,0)\to (\bC^N,0)$ is a formal mapping sending $M$
into $M'$, then either $H(Z)\sim 0$ or the Segre homomorphism
$\Phi_H$ is injective. \endproclaim

\demo{Proof} Let $Z=(z,w)$ and $Z'=(z',w')$ be normal
coordinates for $M$ and $M'$, respectively, as in \S1.2. Thus,
$M$ is given by \thetag{1.2.1} and similarly for $M'$ with $Q$
replaced by $Q'$. We write $H=(F,G)$. We shall prove that if
$H\not\sim 0$, then the mapping $z\mapsto F(z,0)$ is finite,
which implies injectivity of $\Phi_H$ as noted in the proof of
Lemma 6.1. The following is proved in [BR]. Suppose that $M$
and
$M'$ are formal hypersurfaces through $0$ in $\bC^N$ with $M$
essentially finite at $0$. If $H\:(\bC^N,0)\to (\bC^N,0)$ is a
formal mapping sending $M$ into $M'$ and if the transversal
component $G\not\sim 0$ (where $M$ and $M'$ are written in
normal coordinates as above and $H=(F,G)$), then the
formal mapping $z\mapsto F(z,0)$ is finite. Hence, to prove
Proposition 7.1 it suffices to show that $G\sim 0$ implies
$F\sim 0$. 

Thus, assume that $G\sim 0$. The fact that $H$ sends $M$ into
$M'$ is then expressed by 
$$
Q'(F(z,w),\bar F(\chi,\tau),0)\in\Cal I(M).\tag 7.2
$$
We claim that \thetag{7.2} implies 
$$
Q'(F(z,w),\bar
F(\chi,\tau),0)\sim 0.\tag 7.3
$$
To see this, observe that by
linear algebra there exists a nontrivial formal holomorphic
vector field
$$
X=\sum_{j=1}^na_j(z,w)\frac{\partial }{\partial z_j}+
b(z,w)\frac{\partial }{\partial w},\quad a_j,b\in\bC\dbl z,w
\dbr,\tag 7.4
$$
such that
$$ 
XF_j(z,w)\sim 0,\quad j=1,\ldots,n.\tag 7.5
$$
Since $M$ is essentially finite at $0$, $X$ cannot be tangent
to $\Cal I(M)$ (see [BER4, Theorem 11.8.13]; the 
property of having no nontrivial formal holomorphic vector
fields tangent is referred to as holomorphic nondegeneracy of
$M$ at $0$). By \thetag{7.2}, if \thetag{7.3} does not hold,
then there exists
$a\in\bC\dbl z,w,\chi,\tau\dbr$ and an integer $k\geq 1$ such
that 
$$
Q'(F(z,w),\bar
F(\chi,\tau),0)\sim a(z,w,\chi,\tau)\rho(z,w,\chi,\tau)^k,\tag
7.6
$$
where $\rho$ is some generator of $\Cal I(M)$ (e.g.\
$\rho(z,w,\chi,\tau):=w-Q(z,\chi,\tau)$), and $a\not\in\Cal
I(M)$. If we apply
$X$ to \thetag{7.6}, use \thetag{7.5}, and divide by
$\rho^{k-1}$, then we deduce
$$
0\sim (Xa)\rho+ka (X\rho).\tag7.7
$$
We conclude that $a(X\rho)\in\Cal I(M)$ and, hence since $\Cal
I(M)$ is prime, either
$a\in \Cal I(M)$ or $X\rho\in\Cal I(M)$. Both cases are
impossible and consequently the claim \thetag{7.3} is proved.
If
$F(z,w)\not\sim 0$, then there exists a nontrivial formal
curve $s\mapsto (z(s),w(s))$, vanishing at $s=0$, such that
the formal curve
$s\mapsto F(z(s),w(s))$ is nontrivial. Since the formal
curve $s\mapsto (F(z(s),w(s)),0)$ is contained in $M'$ (which
is assumed not to contain a nontrivial formal curve), by
\thetag{7.3}, we conclude that $F(z,w)$ must be identically 0.
This completes the proof of Proposition 7.1.\qed
\enddemo

\demo{Proof of Theorem $2.2$} Since the hypersurface $M$ is
essentially finite at $0$, it is also of
finite type at $0$ (see e.g.\ [BER4, Proposition 9.4.16]). By
the assumption on $M'$ and Lemma 6.3, $M'$ is essentially
finite at
$0$. Moreover, we claim that $M'$ does not contain a nontrivial
formal curve through $0$. Indeed, if the real-analytic
hypersurface
$M'$ were to contain a nontrivial formal curve through $0$,
then
$M'$ would be of infinite type at $0$ in the sense of D'Angelo.
Hence, $M'$ would contain a nontrivial holomorphic subvariety
(see e.g.\ [D'A, Section 3.3.3]), which would contradict the
assumption on $M'$. (One could also apply directly the theorem
of Milman [Mi] which states that if a real-analytic variety
contains a formal variety, then it also contains a holomorphic
one of the same dimension.)  The proof of Theorem 2.2 now
follows from Theorem 2.1 and Proposition 7.1.\qed\enddemo

\demo{Proof of Theorem $1$} Theorem 1 follows
from Theorem 2.2 by using Lemma 6.3.\qed\enddemo

\demo{Proof of Theorem $2.4$}
To prove Theorem 2.4, we denote by $H(Z)$ the formal mapping
$(\bC^N,0)\to (\bC^N,0)$ given by the Taylor series of $h$ at
$0$ (see [BER4, Proposition 1.7.14]). By assumption,
$H(Z)\not\sim 0$. As explained in the proof of Theorem
2.2 above, $M'$ does not contain a
nontrivial formal curve through $0$. We
take normal coordinates for
$M$ and
$M'$ and write
$H=(F,G)$ as in \S1.2. Proposition 7.1
implies that the Segre homomorphism $\Phi_H$ is injective.
Indeed, the proof of Proposition 7.1 in fact shows that 
the mapping $z\mapsto F(z,0)$ is finite. 
It follows  (see e.g.\ [BER4, Theorem 5.1.37]) that
$H$ is not {\it totally degenerate}, i.e.\ 
$$
\det (\partial F_i/\partial z_j)_{1\leq i,j\leq
n}(z,0)\not\sim 0.\tag 7.8
$$
Theorem 2.4 now follows from
the reflection principle in [BR, Theorem 6] (see also [BER4,
Theorem 9.6.1]), since $M'$ is essentially finite at $0$ by
Lemma 6.3.
\qed\enddemo

\subhead Closing remark\endsubhead During the completetion of
this work, the authors became aware of the preprint 
"Convergence of formal 
biholomorphisms between  minimal holomorphically 
nondegenerate real analytic CR
manifolds" (e-print : http://xxx.lanl.gov/ 
abs/math.CV/9901027; 1999) by J. Merker in which
convergence of formal invertible mappings between
holomorphically nondegenerate generic submanifolds of
finite type is claimed. However, the proof in that preprint
has serious flaws. In particular, the proof is based on
Proposition 2.4.6 (in that preprint) which is an incorrect
characterization of holomorphic nondegeneracy for minimal
generic submanifolds. The real algebraic hypersurface in
$\bC^3$ defined by 
$$
\im w=\frac{|z_1|^2|1+z_1\bar z_2|^2}{1+\re(z_1\bar
z_2)}-\re w\frac{\im (z_1\bar z_2)}{1+\re(z_1\bar
z_2)}
$$
gives a counterexample to that proposition.


\Refs\widestnumber\key{BER5}
\ref\key A\by M. Artin\paper On the solution of analytic
equations\jour Invent. Math. \vol 5\yr 1969\pages
277--291\endref

\ref\key AM\by M. F. Atiyah and I. G.
Macdonald\book Introduction to Commutative
Algebra\publ Addison--Wesley\publaddr Reading,
MA\yr 1969\endref

\ref\key BER1\manyby M. S. Baouendi, P. Ebenfelt, and L. P.
Rothschild\paper Algebraicity of holomorphic mappings
between real algebraic sets in $\bC^n$
\jour Acta Math. \vol 177\yr
1996\pages 225--273\endref

\ref\key BER2\bysame\paper CR automorphisms of real
analytic CR manifolds in complex space
\jour Comm. Anal. Geom.\vol 6\pages 291--315\yr 1998\endref

\ref\key BER3\bysame\paper Parametrization of local biholomorphisms of
real analytic hypersurfaces
\jour Asian J. Math.\vol 1\pages 1--16\yr 1997\endref

\ref\key BER4\bysame \book Real Submanifolds
in Complex Space and Their Mappings\publ
Princeton Math. Ser. 47, Princeton
Univ. Press\publaddr Princeton,
NJ\yr 1999\endref
 
\ref\key BER5\bysame\paper Rational dependence of smooth and
analytic  CR mappings 
on their jets\jour Math. Ann.\toappear  \endref 





\ref \key BR\by  M. S. Baouendi and L. P.
Rothschild\paper Geometric properties of mappings between
hypersurfaces in complex space\jour J. Differential Geom.\vol
31\yr 1990\pages 473--499\endref


\ref \key B\by V. K. Beloshapka\paper A uniqueness theorem for
automorphisms of a nondegenerate surface in complex space\jour Math. Notes
47\yr 1990\pages 239--242\endref



\ref\key BK\by E. Brieskorn and H. Kn\"orrer\book Plane
Algebraic Curves\bookinfo(Translated from the German by John
Stillwell)\publ Birkh\"auser Verlag\publaddr Basel-Boston,
Mass.\yr  1986\endref


\ref\key C\by E. Cartan\paper Sur la g\'eom\'etrie
pseudo-conforme des hypersurfaces de deux variables complexes, I\jour
Ann. Math. Pura Appl.\vol 11\yr 1932\pages 17--90\finalinfo (or Oeuvres
II, 1231--1304)\endref

\ref\key D'A\by J. D'Angelo\book Several Complex Variables and
the Geometry of Hypersurfaces\publ Studies in Advanced Math.,
CRC Press\publaddr Boca Raton\yr 1993\endref

\ref\key CM\by S.S. Chern and J. K. Moser\paper Real
hypersurfaces in complex manifolds\jour Acta Math.\vol
133\yr 1974\pages 219--271\endref

%
%
%
%
%
%







\ref\key Go1
\manyby X. Gong
\paper Divergence of the normalization of real-analytic
glancing hypersurfaces
\jour Comm. Partial Differential Eq.\vol 19\yr 1994\pages
643--654
\endref

\ref\key Go2
\bysame
\paper Divergence of the normalization for real Lagrangian
surfaces near complex tangents
\jour Pacific J. Math.\vol 176\yr 1996\pages 311--324
\endref

\ref\key Gu\by R. C. Gunning\book \book Introduction to
Holomorphic Functions of Several Variables.
 \bookinfo Vols. I, II, and III
\publ The Wadsworth \& Brooks/Cole Mathematics Series
\publaddr Pacific Grove, CA \yr 1990\endref

\ref\key Ha \by R. Hartshorne\book Algebraic
Geometry\publ Springer-Verlag\publaddr Berlin\yr
1993\endref
%
%
\ref\key Hu\by X. Huang\paper Schwarz reflection principle in complex spaces of
dimension two\jour  Comm. Partial Differential Eq.\vol 21 
\yr 1996\pages 1781--1828\endref








\ref\key Mel\by R. B. Melrose\paper Equivalence of glancing
hypersurfaces\jour Invent. Math.\vol 37\yr 1976\pages
165--191\endref



\ref\key Mi\by P. Milman\paper Complex analytic and formal
solutions of real analytic equations in $\bC^n$\jour Math.
Ann.\vol 233\yr 1978\pages 1--7\endref 

\ref\key MW\by J. K. Moser and S. M. Webster\paper Normal
forms for real surfaces in $\bC^2$ near complex tangents and
hyperbolic surface transformations\jour Acta Math.\vol 150\yr
1983\pages 255--296\endref

\ref\key O\by T. Oshima\paper On analytic equivalence of
glancing hypersurfaces\jour Sco. Papers College Gen. Ed. Univ.
Tokyo\vol 28\yr 1978\pages 51--57\endref


%



%

\ref\key T
\by  N. Tanaka
\paper On the pseudo-conformal geometry of hypersurfaces of the space of
$n$ complex variables
\jour J. Math. Soc. Japan 
\vol 14 \yr 1962 \pages 397-429
\endref


\ref\key TH\by A. Tumanov and G. M. Henkin
\paper Local characterization of holomorphic
automorphisms of Siegel domains
\jour  Funktsional. Anal. i Prilozhen
\vol 17
\year 1983
\pages 49--61
\transl\nofrills English transl. in \jour Functional Anal.
Appl.
\vol 17\yr 1983\endref

\ref\key VW\by B. L. Van der Waerden\book Modern
Algebra\bookinfo Eighth Edition\publ Springer-Verlag\publaddr
New York, NY\yr 1971\endref



\ref\key W\manyby S. M. Webster\paper Holomorphic symplectic
normalization of a real function\jour Ann. Scuola Norm.
Pisa\vol 19\yr 1992\pages 69--86\endref
%

\ref\key Z\by D. Zaitsev\paper Germs of local automorphisms of
 real-analytic CR structures and analytic dependence on $k$-jets\jour
Math. Research Lett.\vol 4\yr 1997\pages 823--842\endref\smallskip


\ref\key ZS\by O. Zariski and P. Samuel\book Commutative
Algebra\bookinfo Vols. I and II \publ Springer-Verlag\publaddr
New York, NY\yr 1958, 1960\endref 

\endRefs
\enddocument
\end